# RATE OF CONVERGENCE IN THE MULTIDIMENSIONAL CENTRAL LIMIT THEOREM FOR STATIONARY PROCESSES. APPLICATION TO THE KNUDSEN GAS AND TO THE SINAI BILLIARD

By Françoise Pène

*Université de Bretagne Occidentale*

We show how Rio's method [*Probab. Theory Related Fields* **104** (1996) 255–282] can be adapted to establish a rate of convergence in $\frac{1}{\sqrt{n}}$ in the multidimensional central limit theorem for some stationary processes in the sense of the Kantorovich metric. We give two applications of this general result: in the case of the Knudsen gas and in the case of the Sinai billiard.

## 0. Introduction.

0.1. *Context.* We denote probability dynamical system $(\Omega, \mathcal{F}, \nu, T)$ where $(\Omega, \mathcal{F}, \nu)$ is a probability space endowed with a $\nu$-preserving transformation $T$ of $\Omega$. Let a probability dynamical system $(\Omega, \mathcal{F}, \nu, T)$ and a measurable function $f : \Omega \to \mathbb{R}^d$ (with $d \geq 1$) be given; we consider the stationary process $(X_k := f \circ T^k)_{k \geq 1}$.

We say that a $\mathbb{R}^d$-random variable $N$ is Gaussian if, for any $\beta \in \mathbb{R}^d$, the distribution of the real-valued random variable $\langle \beta, N \rangle$ is either a normal distribution or a Dirac measure. With this definition, a Gaussian random variable has a general normal distribution (cf. [16], III-6 for more details). We say that we have a central limit theorem (CLT) for $(X_k)_{k \geq 1}$ if $(\frac{1}{\sqrt{n}} \sum_{k=1}^{n} X_k)_{n \geq 1}$ converges in distribution to a Gaussian random variable $N$.

CLTs in the context of dynamical systems have been established in many articles (cf. [13, 27, 36, 39]). In these works, multidimensional central limit theorem follows directly from one-dimensional central limit theorem. Indeed, let us recall that the fact that $(\frac{1}{\sqrt{n}} \sum_{k=1}^{n} X_k)_{n \geq 1}$ converges in dis-tribution to $N$ means that, for all $\beta \in \mathbb{R}^d$, the one-dimensional process $(\frac{1}{\sqrt{n}} \sum_{k=1}^n \langle \beta, X_k \rangle)_{n \geq 1}$ converges in distribution to $\langle \beta, N \rangle$.

In the present paper, we are interested in questions of speed of convergence in the CLT for multidimensional stationary processes. There are many ways of estimating the speed in the CLT. Let us endow $\mathbb{R}^d$ with the supremum norm $|\cdot|_\infty$ and with the Borel $\sigma$-algebra $\mathcal{B}(\mathbb{R}^d)$. For any $\phi \colon \mathbb{R}^d \to \mathbb{R}$, we denote by $L_\phi$ its Lipschitz constant:

$$L_\phi := \sup_{x,y \in \mathbb{R}^d \,:\, x \neq y} \frac{|\phi(x) - \phi(y)|}{|x - y|_\infty}.$$

We can estimate the rate of convergence in the CLT for $(X_k)_k$ by estimating the following quantities:

(a) uniform norm of the difference between the distribution functions ($DF$ metric):

$$DF_n := \sup_{(x_1,\ldots,x_d) \in \mathbb{R}^d} \left| \nu\left( \frac{1}{\sqrt{n}} \sum_{k=1}^n X_k^{(1)} \leq x_1, \ldots, \frac{1}{\sqrt{n}} \sum_{k=1}^n X_k^{(d)} \leq x_d \right) \right.$$
$$\left. - \mathbb{P}(N^{(1)} \leq x_1, \ldots, N^{(d)} \leq x_d) \right|,$$

where $X_k^{(i)}$ and $N^{(i)}$ are the $i$th coordinates of $X_k$ and of $N$, respectively;

(b) in the sense of the Prokhorov metric ($\Pi$ metric):

$$\Pi_n := \inf\left\{ \varepsilon > 0 \colon \forall\, B \in \mathcal{B}(\mathbb{R}^d), \nu\left( \frac{1}{\sqrt{n}} \sum_{k=1}^n X_k \in B \right) - \mathbb{P}(N \in B^\varepsilon) \leq \varepsilon \right\},$$

where $B^\varepsilon$ is the open $\varepsilon$-neighborhood of $B$;

(c) in the sense of the Lipschitz bounded metric ($LB$ metric):

$$LB_n := \sup\left\{ \left| \mathbb{E}_\nu\left[ \phi\left( \frac{1}{\sqrt{n}} \sum_{k=1}^n X_k \right) \right] - \mathbb{E}[\phi(N)] \right|, \phi \colon \mathbb{R}^d \to \mathbb{R}, \|\phi\|_\infty + L_\phi \leq 1 \right\};$$

(d) in the sense of the Kantorovich metric ($\kappa$ metric):

$$\kappa_n := \sup\left\{ \left| \mathbb{E}_\nu\left[ \phi\left( \frac{1}{\sqrt{n}} \sum_{k=1}^n X_k \right) \right] - \mathbb{E}[\phi(N)] \right|, \phi \colon \mathbb{R}^d \to \mathbb{R}, L_\phi \leq 1 \right\}.$$

We will give additional details about the metrics corresponding to these quantities.

### 0.2. *Previous results.*

0.2.1. *One-dimensional processes* $(d = 1)$. When $d = 1$, it is classical to estimate the speed in the CLT in the sense of the uniform error between the distribution functions ($DF$ metric).

In [1, 3, 15] a rate of convergence in $\frac{1}{\sqrt{n}}$ in the sense of the $DF$ metric is established for sequences $(X_k)_k$ of independent identically distributed random variables such that $\mathbb{E}[|X_1|^3] < +\infty$. Moreover, this result is optimal.

This result has been extended to some martingale processes (cf. [7]) and to some stationary processes ([37] extended in [22, 23] and in [25, 26]).

When $d = 1$ and when $(X_k)_k$ is a sequence of independent identically distributed random variables, Nagaev [31] establishes a nonuniform estimate for the difference between the distribution functions:

$$\exists L > 0 \ \forall n \geq 1 \qquad \left| \nu\left( \frac{X_1 + \cdots + X_n}{\sqrt{n}} \leq x \right) - \mathbb{P}(N \leq x) \right| \leq \frac{L\mathbb{E}[|X_1|^3]}{\sqrt{n}(1 + |x|^3)}.$$

A direct consequence of this is a speed of convergence in $\frac{1}{\sqrt{n}}$ in the sense of the $DF$ metric but also in the sense of the Kantorovich metric (cf. Proposition 0.10). Results in the sense of the Kantorovich metric have been established by many authors. Let us mention the article [41] of Sunklodas.

0.2.2. *Multidimensional processes.* For sequences of independent identically distributed random variables $(X_k)_{k \geq 1}$ with values in $\mathbb{R}^d$ and admitting an invertible covariance matrix and admitting moments of the third order, a speed in $\frac{1}{\sqrt{n}}$ is established by Bergström in the sense of the $DF$ metric (cf. theorem of page 121 in [2]). This result gives an extension of the Berry–Esseen result to the $d$-dimensional case.

In [23] Jan shows that Rio's result can be extended to the $d$-dimensional case (in the sense of the $DF$ metric).

In [44] Yurinskii establishes an inequality linking the Prokhorov metric with characteristic functions. This result allows to establish a rate of convergence in $\frac{1}{\sqrt{n}}$ for sequences of independent identically distributed random variables in the sense of the Prokhorov metric. We will recall and use this inequality in the case of the Knudsen gas.

Let us also mention the works of [4, 35, 38, 42] in which the rate of convergence in the CLT is estimated in other ways for independent random variables sequences.

0.3. *Contents of the present paper.* In Section 1 we give a result of speed of convergence in the CLT for some stationary processes $(X_k)_k$ in the sense of the Kantorovich metric (Theorem 1.1). The proof of this result, given in the Appendix, uses an adaptation of the method developed in [37] and extended in [25, 26]. In these papers, a rate of convergence in the sense of the $DF$ metric has been established in the one-dimensional case. Our hypothesis is analogous to the hypothesis of [25, 26] but weaker than it.

In Sections 2 and 3 we will give applications of our result.

In Section 2 we study the Knudsen gas model studied by Boatto and Golse in [6]. We use a Markov model of it. We show that the results of [6] are related to questions of rate of convergence in the CLT. For this model, we will, first, estimate the speed of convergence in the CLT in the sense of the Prokhorov metric (using Yurinskii's result of [44] and an extension of the method of perturbation of quasi-compact operators [20, 29, 30]). Second, we apply Theorem 1.1 of Section 1 and establish a rate in the sense of the Kantorovich metric. This result gives an extension of a result of [6].

In Section 3 we are interested in the question of the rate of convergence in the multidimensional CLT for $(X_k := f \circ T^k)_k$ where $T$ is the billiard transformation of the Sinai billiard [40] and $f$ is a smooth (Hölder continuous) function. Using Theorem 1.1 of Section 1, we establish a rate of convergence in $\frac{1}{\sqrt{n}}$ in the sense of the Kantorovich metric. This is, to our knowledge, the first time that a speed of convergence in $\frac{1}{\sqrt{n}}$ is established in this context. In [33, 34] a speed of convergence in $\frac{1}{n^{1/2-\alpha}}$ for any $\alpha > 0$ has been established in the sense of the $DF$ metric and in the sense of the Prokhorov metric (with an adaptation of a method developed by Jan in [23] using characteristic functions). The result we present here does not exactly improve [33, 34]. It gives a better rate with another metric.

0.4. *Some metrics for probability measures on $\mathbb{R}^d$.* Let us denote by $\mathcal{M}_1(\mathbb{R}^d)$ the set of probability measures on $(\mathbb{R}^d, \mathcal{B}(\mathbb{R}^d))$, where $\mathcal{B}(\mathbb{R}^d)$ is the Borel $\sigma$-algebra of $\mathbb{R}^d$.

DEFINITION 0.1 (The $DF$ metric).   For all $P, Q$ in $\mathcal{M}_1(\mathbb{R}^d)$, we define

$$DF(P,Q) := \sup_{(x_1,\dots,x_d) \in \mathbb{R}^d} \left| P\left( \prod_{i=1}^{d} ]-\infty; x_i] \right) - Q\left( \prod_{i=1}^{d} ]-\infty; x_i] \right) \right|.$$

DEFINITION 0.2 (The Prokhorov metric; cf. [5, 14]).   For all $P, Q$ in $\mathcal{M}_1(\mathbb{R}^d)$, we define

$$\Pi(P,Q) := \inf\{\varepsilon > 0 : \forall B \in \mathcal{B}(\mathbb{R}^d), (P(B) - Q(B^\varepsilon)) \leq \varepsilon\}.$$

DEFINITION 0.3 (The Ky Fan metric for random variables).   If $X$ and $Y$ are two $\mathbb{R}^d$-valued random variables defined on the same probability space $(E_0, \mathcal{T}_0, \mathbb{P}_0)$, we define

$$\mathcal{K}(X,Y) := \inf\{\varepsilon > 0 : \mathbb{P}_0(|X - Y|_\infty > \varepsilon) < \varepsilon\}.$$

Let us recall that $\lim_{n\to+\infty}\mathcal{K}(X_n,Y)=0$ means that $(X_n)_n$ converges in probability to $Y$.

PROPOSITION 0.4 ([14], Corollary 11.6.4). *For all $P,Q$ in $\mathcal{M}_1(\mathbb{R}^d)$, the quantity $\Pi(P,Q)$ is the infimum of $\mathcal{K}(X,Y)$ where $X$ and $Y$ are two $\mathbb{R}^d$-valued random variables defined on the same probability space and such that the distribution of $X$ is $P$ and the distribution of $Y$ is $Q$.*

DEFINITION 0.5 (The bounded Lipschitz metric). For all $P,Q$ in $\mathcal{M}_1(\mathbb{R}^d)$, we define

$$BL(P,Q):=\sup\{|\mathbb{E}_P[\phi]-\mathbb{E}_Q[\phi]|,\phi\colon\mathbb{R}^d\to\mathbb{R},\|\phi\|_\infty+L_\phi\leq 1\}.$$

In particular, for any bounded Lipschitz continuous function $\phi\colon\mathbb{R}^d\to\mathbb{R}$, we have

$$|\mathbb{E}_P[\phi]-\mathbb{E}_Q[\phi]|\leq BL(P,Q)\times(\|\phi\|_\infty+L_\phi).$$

PROPOSITION 0.6 ([14], Theorem 11.3.3). *Let $(P_n)_n$ be a sequence of $\mathcal{M}_1(\mathbb{R}^d)$ and let $P$ be in $\mathcal{M}_1(\mathbb{R}^d)$. The following properties are equivalent:*

(i) *the sequence $(P_n)_n$ of probability measures converges weakly to $P$;*
(ii) $\lim_{n\to+\infty}\Pi(P_n,P)=0$;
(iii) $\lim_{n\to+\infty}BL(P_n,P)=0$;
(iv) $\lim_{n\to+\infty}DF(P_n,P)=0$, *if $P$ has a continuous distribution.*

More precisely, we have (cf. [28], Proposition 1.2 and [14], Problem 11.3.5):

$$\tfrac{1}{3}BL(P,Q)\leq\Pi(P,Q)\leq(\tfrac{3}{2}BL(P,Q))^{1/3}.$$

Let us denote by $\mathcal{M}_{1,\mathrm{int}}(\mathbb{R}^d)$ the set of probability measures on $(\mathbb{R}^d,\mathcal{B}(\mathbb{R}^d))$ admitting moments of the first order.

DEFINITION 0.7 (The Kantorovich metric, cf. [14, 16]). For all $P,Q$ in $\mathcal{M}_{1,\mathrm{int}}(\mathbb{R}^d)$, we define

$$\kappa(P,Q):=\sup\{|\mathbb{E}_P[\phi]-\mathbb{E}_Q[\phi]|,\phi\colon\mathbb{R}^d\to\mathbb{R},L_\phi\leq 1\}.$$

In particular, for any Lipschitz continuous function $\phi\colon\mathbb{R}^d\to\mathbb{R}$, we have

$$|\mathbb{E}_P[\phi]-\mathbb{E}_Q[\phi]|\leq\kappa(P,Q)\times L_\phi.$$

PROPOSITION 0.8 ([14], Theorem 11.8.2). *For all $P, Q$ in $\mathcal{M}_1(\mathbb{R}^d)$, the quantity $\kappa(P, Q)$ is the infimum of $\mathbb{E}[|X - Y|_\infty]$, where $X$ and $Y$ are two $\mathbb{R}^d$-valued random variables defined on the same probability space and such that the distribution of $X$ is $P$ and the distribution of $Y$ is $Q$.*

PROPOSITION 0.9. *Let $(P_n)_n$ be a sequence of $\mathcal{M}_{1,\mathrm{int}}(\mathbb{R}^d)$ and $P$ in $\mathcal{M}_{1,\mathrm{int}}(\mathbb{R}^d)$. The following properties are equivalent:*

(i) *the sequence $(P_n)_n$ of probability measures converges weakly to $P$ and we have $\lim_{n \to +\infty} \int_{\mathbb{R}^d} |x|_\infty \, dP_n(x) = \int_{\mathbb{R}^d} |x|_\infty \, dP(x)$;*

(ii) *$\lim_{n \to +\infty} \kappa(P_n, P) = 0$.*

PROOF. According to Proposition 0.6 and to the fact that $BL(P_n, P) \le \kappa(P_n, P)$, it is easy to see that (ii) implies (i).

Let us now suppose that (i) is true and prove that (ii) is then true. Let us write $\alpha_n := |\int_{\mathbb{R}^d} |x|_\infty \, dP_n(x) - \int_{\mathbb{R}^d} |x|_\infty \, dP(x)|$. Let $\phi : \mathbb{R}^d \to \mathbb{R}$ be any Lipschitz continuous function with Lipschitz constant $L_\phi$ bounded by 1. For any nonnegative real number $M$, we define $\psi_M : \mathbb{R}^d \to \mathbb{R}$ by

$$\psi_M(x) = \begin{cases} \phi(x), & \text{if } |\phi(x) - \phi(0)| \le M, \\ \phi(0) + M, & \text{if } \phi(x) \ge \phi(0) + M, \\ \phi(0) - M, & \text{if } \phi(x) \le \phi(0) - M. \end{cases}$$

For all $M > 0$ and all integer $n \ge 0$, we have

$$|\mathbb{E}_{P_n}[\psi_M] - \mathbb{E}_P[\psi_M]| \le (M+1) BL(P_n, P)$$

and

$$\forall x \in \mathbb{R}^d \qquad |\psi_M(x) - \phi(x)| = |\psi_M(x) - \phi(x)| \mathbf{1}_{\{|x|_\infty > M\}}$$
$$\le (|x|_\infty - M) \mathbf{1}_{\{|x|_\infty > M\}}$$

and therefore

$$|\mathbb{E}_{P_n}[\psi_M - \phi]| \le \mathbb{E}_{P_n}[(|\cdot|_\infty - M) \mathbf{1}_{\{|\cdot|_\infty > M\}}]$$
$$\le \mathbb{E}_{P_n}[|\cdot|_\infty - \min(|\cdot|_\infty, M)]$$
$$\le \mathbb{E}_P[|\cdot|_\infty] + \alpha_n - \mathbb{E}_P[\min(|\cdot|_\infty, M)] + (M+1) BL(P_n, P)$$
$$\le \mathbb{E}_P[(|\cdot|_\infty - M) \mathbf{1}_{\{|\cdot|_\infty > M\}}] + \alpha_n + (M+1) BL(P_n, P).$$

Hence, for any $M > 0$ and any integer $n \ge 0$, we have

$$|\mathbb{E}_{P_n}[\phi] - \mathbb{E}_P[\phi]| \le 2\mathbb{E}_P[(|\cdot|_\infty - M) \mathbf{1}_{\{|\cdot|_\infty > M\}}] + \alpha_n + 2(M+1) BL(P_n, P).$$

Let a real number $\varepsilon > 0$ be given. Let us fix $M_\varepsilon$ such that $\mathbb{E}_P[|\cdot|_\infty \mathbf{1}_{\{|\cdot|_\infty > M_\varepsilon\}}] < \frac{\varepsilon}{5}$. Let $N_\varepsilon$ be such that, for any integer $n \ge N_\varepsilon$, we have $\alpha_n \le \frac{\varepsilon}{5}$ and $(M_\varepsilon +$

1)$BL(P_n, P) \leq \frac{\varepsilon}{5}$ (according to Proposition 0.6, such a $N_\varepsilon$ exists). Therefore, for any $n \geq N_\varepsilon$ and any Lipschitz continuous function $\phi \colon \mathbb{R}^d \to \mathbb{R}$ with Lipschitz constant $L_\phi$ less than 1, we have

$$|\mathbb{E}_{P_n}[\phi] - \mathbb{E}_P[\phi]| < \varepsilon. \qquad \square$$

PROPOSITION 0.10 (The Kantorovich metric, case $d = 1$ ([14], problem 11.8.1)). *For all $P, Q$ in $\mathcal{M}_{1,\mathrm{int}}(\mathbb{R})$, we have*

$$\kappa(P, Q) = \int_{\mathbb{R}} |P(]-\infty; x]) - Q(]-\infty; x])| \, dx.$$

**1. Abstract theorem.** We will denote by $\mathbf{Lip}(\mathbb{R}^d, \mathbb{R})$ the set of Lipschitz continuous functions from $\mathbb{R}^d$ into $\mathbb{R}$. Let a probability space $(\Omega, \mathcal{F}, \nu)$ be given. For any $\nu$-integrable function $f \colon \Omega \to \mathbb{R}$, we denote by $\mathbb{E}_\nu[f]$ the expectation of $f$ with respect to probability measure $\nu$:

$$\mathbb{E}_\nu[f] := \int_\Omega f \, d\nu.$$

For all real-valued functions $f, g$ in $L^2(\Omega, \nu)$, we recall the definition of the covariance of $f$ and $g$ (with respect to measure $\nu$):

$$\mathrm{Cov}_\nu(f, g) = \mathbb{E}_\nu[fg] - \mathbb{E}_\nu[f]\mathbb{E}_\nu[g].$$

For all $A = (a_1, \ldots, a_d)$ and $B = (b_1, \ldots, b_d)$ in $\mathbb{R}^d$, we denote by $A \otimes B$ and $A^{\otimes 2}$ the square matrices given by

$$A \otimes B := (a_i b_j)_{i,j=1,\ldots,d} \quad \text{and} \quad A^{\otimes 2} := A \otimes A.$$

Let $M = (M_{i,j})_{i,j=1,\ldots,d}$ be a random variable on $\Omega$ with values in the set of the square matrices such that, for any $i, j = 1, \ldots, d$, $M_{i,j}$ is $\nu$-integrable. Then, the expectation $\mathbb{E}_\nu[M]$ of $M$ is the $d$-dimensional matrix given by

$$\mathbb{E}_\nu[M] := (\mathbb{E}_\nu[M_{i,j}])_{i,j=1,\ldots,d}.$$

Let us consider a sequence of stationary $\mathbb{R}^d$-valued random variables $(X_k)_{k \geq 0}$ defined on $(\Omega, \mathcal{F}, \nu)$.

For any integer $n \geq 1$, we write $S_n := \sum_{k=1}^n X_k$ and $S_0 = 0$.

Using an adaptation of Rio's method developed in [37], we will establish the following result:

THEOREM 1.1. *Let $(X_k)_{k \geq 0}$ be a sequence of stationary $\mathbb{R}^d$-valued bounded random variables defined on $(\Omega, \mathcal{F}, \nu)$ with expectation 0. Let us suppose that there exist two real numbers $C \geq 1$, $M \geq \max(1, \|X_0\|_\infty)$ and an integer $r \geq 0$ and a sequence of real numbers $(\varphi_{p,l})_{p,l}$ bounded by 1 with $\sum_{p \geq 1} p \max_{l=0,\ldots,\lfloor p/(r+1)\rfloor} \varphi_{p,l} < +\infty$ such that for any integers $a, b, c \geq 0$ satisfying $1 \leq a + b + c \leq 3$, for any integers $i, j, k, p, q, l$ with $1 \leq i \leq j \leq$* $k \le k + p \le k + p + q \le k + p + l$, for any $i_1, i_2, i_3 \in \{1, \dots, d\}$, for any $F : \mathbb{R}^d \times ([-M; M]^d)^3 \to \mathbb{R}$ bounded, differentiable, with bounded differential, we have

$$
\begin{aligned}
(1) \quad & |\mathrm{Cov}(F(S_{i-1}, X_i, X_j, X_k), (X_{k+p}^{(i_1)})^a (X_{k+p+q}^{(i_2)})^b (X_{k+p+l}^{(i_3)})^c)| \\
& \le C(\|F\|_{L^\infty} + \||DF|_\infty\|_{L^\infty}) \varphi_{p,l},
\end{aligned}
$$

where $DF$ is the Jacobian matrix of $F$ and $X_m^{(s)}$ is the sth coordinate of $X_m$. Then, the following limit exists:

$$
\Sigma^2 := \lim_{n \to +\infty} \frac{1}{n} (\mathbb{E}[S_n^{\otimes 2}]).
$$

If $\Sigma^2 = 0$, then the sequence $(S_n)_n$ is bounded in $L^2$.

Otherwise the sequence of random variables $(\frac{S_n}{\sqrt{n}})_{n \ge 1}$ converges in distribution to a Gaussian random variable $N$ with expectation 0 and with covariance matrix $\Sigma^2$ and there exists a real number $B > 0$ such that, for any integer $n \ge 1$ and any Lipschitz continuous function $\phi : \mathbb{R}^d \to \mathbb{R}$, we have

$$
\left| \mathbb{E}\left[ \phi\left( \frac{S_n}{\sqrt{n}} \right) \right] - \mathbb{E}[\phi(N)] \right| \le \frac{B L_\phi}{\sqrt{n}}.
$$

The proof of this theorem is given in the Appendix.

## 2. Application to the Knudsen gas.
Following Boatto and Golse [6], we are interested in a generalized model of the Knudsen gas with an isotropic component. In the present section we use a probabilistic approach. We show how this problem can be modeled by a Markov chain. Using the method of perturbation of operators due to Nagaev (see [19, 20, 29, 30]), we get a rate in $\frac{1}{\sqrt{n}}$ in the multidimensional CLT in the sense of the Prokhorov metric. Moreover, we establish the same rate for the Kantorovich metric. This second result is an application of Theorem 1.1 and gives an extension of a theorem of [6].

### 2.1. The model.
In this section we will make the following assumption.

HYPOTHESIS 2.1. $(\Omega, \mathcal{F}, \nu, T)$ is an invertible probability dynamical system, $a : \Omega \to \mathbb{R}^d$ is a $\nu$-centered square integrable function and $\alpha$ is a fixed real number satisfying $0 < \alpha < 1$.

The invertibility hypothesis is not restrictive since any dynamical system admits an invertible extension (its natural extension). Moreover, let us recall that any stationary sequence of centered and square integrable random variables admits a representation of the form $(Y_k = a \circ T^k)_k$ with $(\Omega, \mathcal{F}, \nu, T)$ and

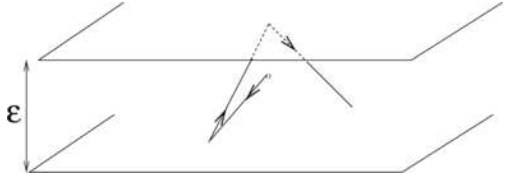

<small>Fig. 1.</small>

$a$ as in Hypothesis 2.1. We denote by $L^p(\Omega, \mathbb{R}^d)$ the set of measurable functions $f \colon \Omega \to \mathbb{R}^d$ such that $\int_\Omega |f(\omega)|_\infty^p \, d\nu(\omega) < +\infty$. For any $f \in L^p(\Omega, \mathbb{R}^d)$, we define $\|f\|_{L^p} := (\int_\Omega |f(\omega)|_\infty^p \, d\nu(\omega))^{1/p}$. We denote by $L^\infty(\Omega, \mathbb{R}^d)$ the set of measurable functions $f \colon \Omega \to \mathbb{R}^d$ which are $\nu$-almost surely bounded by some constant and, for such a function, we denote by $\|f\|_\infty$ the following real number:

$$\|f\|_\infty := \inf\{M > 0 \colon \nu(\{\omega \in \Omega \colon |f(\omega)|_\infty > M\}) = 0\}.$$

We consider a system of particles moving independently in $\mathbb{R}^{d+1}$ between two $d$-dimensional horizontal plates $\mathbb{R}^d \times \{0\}$ and $\mathbb{R}^d \times \{\varepsilon\}$ separated by some small distance $\varepsilon > 0$. We suppose that these particles move with speed $\frac{1}{\varepsilon}(a(\omega), \pm 1)$ parametrized by $\omega$. In our model, the speed and the parameter $\omega$ only change when the particle hits one of the plates; a particle incoming to the upper plate with the speed $\frac{1}{\varepsilon}(a(\omega), 1)$ will outgo:

(a) either (with probability $1 - \alpha$) with the speed $\frac{1}{\varepsilon}(a(T(\omega)), -1)$;

(b) or (with probability $\alpha$) with the speed $\frac{1}{\varepsilon}(a(\omega'), -1)$, where $\omega'$ is given by a random variable (independent of $\omega$) with distribution $\nu$.

We make analogous assumptions for reflections off the lower plate (replacing $\frac{1}{\varepsilon}$ by $-\frac{1}{\varepsilon}$ and $-\frac{1}{\varepsilon}$ by $\frac{1}{\varepsilon}$). We are interested in the behavior of such a model when $\varepsilon$ goes to zero. See Figure 1.

Let us study the evolution of a single particle moving in this system. Let us write $\delta = 1$ if, at time 0, the particle is pointing upward and $\delta = -1$ if, at time 0, the particle is pointing downward. Then, the speed of the particle between the $n$th and the $(n+1)$st collision off one of the plates is $\frac{1}{\varepsilon}(a(X_n), \delta(-1)^n)$, where $(X_n)_n$ is a Markov chain such that the conditional law of $X_{n+1}$ with respect to $(X_0, \dots, X_n)$ is $(1-\alpha)\delta_{T(X_n)} + \alpha\nu$. Let us notice that the measure $\nu$ is an invariant probability measure for this Markov chain. More precisely, we define $(X_n)_{n \in \mathbb{Z}}$ as follows.

NOTATION 2.1. We consider the probability space $(\tilde{\Omega}, \tilde{\mathcal{F}}, \tilde{\nu})$ with $\tilde{\Omega} := \Omega^{\mathbb{Z}}$, $\tilde{\mathcal{F}}$ the product $\sigma$-algebra and $\tilde{\nu}$ the unique probability measure defined on $\tilde{\Omega}$ such that we have

$$\mathbb{E}_{\tilde{\nu}}[f(X_n)] = \int_\Omega f \, d\nu$$

and

$$\mathbb{E}_{\tilde{\nu}}[f(X_{n+1})|X_n, X_{n-1}, \dots] = (1 - \alpha)f(T(X_n)) + \alpha \int_{\Omega} f \, d\nu,$$

with $X_n : \tilde{\Omega} \to \Omega$ given by $X_n((\omega_m)_{m \in \mathbb{Z}}) := \omega_n$. We define the transformation $\tilde{T}$ on $\tilde{\Omega}$ by $\tilde{T}((\omega_n)_{n \in \mathbb{Z}}) = (\omega_{n+1})_{n \in \mathbb{Z}}$.

Since $\nu$ is $T$-invariant, the existence of $\tilde{\nu}$ is a consequence of a result of Ionescu Tulcea (cf. [21], [32], page 154). With this notation, if the particle is at time 0 at the position $(x, \varepsilon z)$ (with $x \in \mathbb{R}^d$ and $z \in [0; 1]$) with the speed $\frac{1}{\varepsilon}(a(X_0), -1)$ parametrized by $X_0$, then its horizontal position at time $s > 0$ will be given by

$$\xi_{\varepsilon}^-(s, x, z, \cdot) := x + \varepsilon \cdot z \cdot a(X_0) + \varepsilon \sum_{k=1}^{\lfloor (s/\varepsilon^2) - z \rfloor} a(X_k) + \varepsilon \cdot \left\{ \frac{s}{\varepsilon^2} - z \right\} a(X_{\lfloor (s/\varepsilon^2) - z \rfloor + 1}),$$

where $\{u\}$ is the fractional part of $u$. For symmetry reasons, if the particle was at time 0 at the position $(x, \varepsilon z)$ (with $x \in \mathbb{R}^d$ and $z \in [0; 1]$) with the speed $\frac{1}{\varepsilon}(a(X_0), 1)$ parametrized by $X_0$, then its horizontal position at time $s > 0$ is given by $\xi_{\varepsilon}^+(s, x, z, \cdot) := \xi_{\varepsilon}^-(s, x, 1 - z, \cdot)$.

2.2. *Results.* In [6] Boatto and Golse have studied the following quantities:

$$F_{\varepsilon, \phi}^{\pm}(s, x, z, \omega) = \mathbb{E}_{\tilde{\nu}}[\phi(\xi_{\varepsilon}^{\pm}(s, x, z, \cdot))|X_0 = \omega]$$

(their $f_{\alpha, \varepsilon}^{\pm}$ corresponds to our $F_{\varepsilon, \psi}^{\mp}$ with $-a$ instead of $a$). More precisely, they establish the following result in the situation when the dynamical system is given by the algebraic automorphism $T_0$ of the two-dimensional torus $\mathbb{T}^2 = \frac{\mathbb{R}^2}{\mathbb{Z}^2}$ given by the matrix $\begin{pmatrix} 2 & 1 \\ 1 & 1 \end{pmatrix}$. We recall that $T_0$ preserves the Haar measure $\nu_0$ on $\mathbb{T}^2$.

THEOREM 2.2.1 ([6]). *Let us suppose that* $(\Omega, \mathcal{F}, \nu, T) = (\mathbb{T}^2, \mathcal{B}(\mathbb{T}^2), \nu_0, T_0)$. *Let* $\phi$ *be a smooth bounded function with bounded derivatives up to order 4. Let* $a$ *be a* $\nu$-*centered function belonging to* $H^{\beta}(\mathbb{T}^2, \mathbb{R}^d)$ *with* $\beta > 1$ *such that the matrix* $D(a) := \sum_{k \in \mathbb{Z}} (1 - \alpha)^{|k|} \mathbb{E}_{\nu}[a \otimes a \circ T^k]$ *is invertible. Then, for any real number* $t_0 > 0$, *we have*

$$\sup_{s \in [0; t_0]} \sup_{x \in \mathbb{R}^d} \sup_{\omega \in \Omega} \sup_{z \in [0; 1]} |F_{\varepsilon, \phi}^{\pm}(s, x, z, \omega) - \mathbb{E}[\phi(x + B_s)]| = O(\varepsilon),$$

*where* $(B_s)_{s \in \mathbb{R}}$ *is a* $d$-*dimensional Brownian motion with zero mean and with covariance matrix* $D(a)$.

We will show how this result is related to the central limit theorem for $(a(X_k))_k$ and give some extensions of it. Indeed, for any real number $\beta > 1$, the Sobolev space $H^\beta(\mathbb{T}^2, \mathbb{R}^d)$ is contained in $L^\infty(\mathbb{T}^2, \mathbb{R}^d)$ and we have

REMARK 2.2.2. *Under Hypothesis 2.1, if $a$ is in $L^\infty(\Omega, \mathbb{R}^d)$ and $\phi : \mathbb{R}^d \to \mathbb{R}$ is a Lipschitz continuous function, then we have*

$$\sup_{s>0, x \in \mathbb{R}^d, z \in [0;1]} \left\| F_{\varepsilon,\phi}^-(s,x,z,\omega) - \mathbb{E}_{\tilde\nu}\left[ \phi\left( x + \varepsilon \sum_{k=0}^{\lfloor s/\varepsilon^2 \rfloor} a(X_k) \right) \right] \right\|_{L^\infty(\Omega)}$$

$$\leq \varepsilon L_\phi \|a\|_\infty \left( 4 + 2 \sum_{l \geq 1} l(1-\alpha)^{l-1} \right).$$

PROOF. We have

$$\left| F_{\varepsilon,\phi}^-(s,x,z,\omega) - F_{\varepsilon,\phi}^-\left( \left( \left\lfloor \frac{s}{\varepsilon^2} \right\rfloor + 1 \right) \varepsilon^2, x, 1, \omega \right) \right| \leq 4\varepsilon L_\phi \|a\|_{L^\infty}.$$

For any $k \geq 1$, we have

$$F_{\varepsilon,\phi}^-((k+1)\varepsilon^2, x, 1, \omega) = \sum_{j=0}^k \alpha^j (1-\alpha)^{k-j} \sum_{l_0 \geq 0; l_1 \geq 1, \ldots, l_j \geq 1 \,:\, l_0 + \cdots + l_j = k} \alpha_{l_0, \ldots, l_j}(\omega),$$

with

$$\alpha_k(\omega) := \phi\left( x + \varepsilon \sum_{m=0}^k a(T^m(\omega)) \right)$$

and

$$\alpha_{l_0, \ldots, l_j}(\omega) := \int_\Omega \cdots \int_\Omega \phi\left( x + \varepsilon \left[ \left( \sum_{m=0}^{l_0} a(T^m(\omega)) \right) \right. \right.$$
$$\left. \left. + \sum_{i=1}^j \left( \sum_{m_i=1}^{l_i} a(T^{m_i}(\omega_i)) \right) \right] \right) d\nu(\omega_1) \cdots d\nu(\omega_j),$$

from which we deduce that we have

$$|F_{\varepsilon,\varphi}^-((k+1)\varepsilon^2, x, 1, \omega) - F_{\varepsilon,\varphi}^-((k+1)\varepsilon^2, x, 1, \omega')|$$

$$\leq 2\varepsilon L_\phi \|a\|_\infty \sum_{l_0 \geq 0} (l_0 + 1)(1-\alpha)^{l_0}.$$

Moreover, we have

$$\mathbb{E}_{\tilde\nu}\left[ F_{\varepsilon,\phi}^-\left( \left( \left\lfloor \frac{s}{\varepsilon^2} \right\rfloor + 1 \right) \varepsilon^2, x, 1, \cdot \right) \right] = \mathbb{E}_{\tilde\nu}\left[ \phi\left( x + \varepsilon \sum_{k=0}^{\lfloor s/\varepsilon^2 \rfloor} a(X_k) \right) \right]. \qquad \square$$

PROPOSITION 2.2.3.   *Let us suppose that Hypothesis* 2.1 *is satisfied. For any integer* $k \geq 0$, *we have* $\mathbb{E}_{\tilde{\nu}}[a(X_0) \otimes a(X_k)] = (1-\alpha)^k \mathbb{E}_{\nu}[a \otimes a \circ T^k]$. *Moreover, the following limit exists:*

$$D(a) := \lim_{n \to +\infty} \mathbb{E}_{\tilde{\nu}}\left[\left(\frac{1}{\sqrt{n}}\sum_{k=0}^{n-1} X_k\right)^{\otimes 2}\right]$$

*and satisfies*

$$D(a) = \sum_{k \in \mathbb{Z}} (1-\alpha)^{|k|} \mathbb{E}_{\nu}[a \otimes a \circ T^k].$$

PROOF.   This is a consequence of the fact that we have $\mathbb{E}_{\tilde{\nu}}[a(X_k)|X_0] = (1-\alpha)^k a(T^k(X_0))$.   □

Here, we prove the two following results:

THEOREM 2.2.4 (Rate of convergence in the CLT in the sense of Prokhorov). *Under Hypothesis* 2.1, *if* $a : \Omega \to \mathbb{R}^d$ *belongs to* $L^3(\Omega, \mathbb{R}^d)$ *and to* $L^{\lfloor d/2 \rfloor + 1}$, *then* $(\frac{1}{\sqrt{n}}\sum_{k=0}^{n-1} a(X_k))_n$ *converges in distribution to a centered Gaussian random variable with covariance matrix* $D(a) = \sum_{k \in \mathbb{Z}}(1-\alpha)^{|k|}\mathbb{E}_{\nu}[a \otimes a \circ T^k]$. *Moreover there exists a real number* $A > 0$ *such that, for any integer* $n \geq 1$, *we have*

$$\Pi\left(\tilde{\nu}_*\left(\frac{1}{\sqrt{n}}\sum_{k=0}^{n-1} a(X_k)\right), \mathcal{N}(0, D(a))\right) \leq \frac{A}{\sqrt{n}},$$

*where* $\tilde{\nu}_*(\frac{1}{\sqrt{n}}\sum_{k=0}^{n-1} a(X_k))$ *denotes the distribution of* $\frac{1}{\sqrt{n}}\sum_{k=0}^{n-1} a(X_k)$ *with respect to* $\tilde{\nu}$.

THEOREM 2.2.5 (Rate of convergence in the CLT in the sense of the Kantorovich metric).   *Under Hypothesis* 2.1, *if* $a : \Omega \to \mathbb{R}^d$ *is a* $\nu$-*centered function belonging to* $L^{\infty}(\Omega, \mathbb{R}^d)$, *then there exists a constant* $B > 0$ *such that, for any Lipschitz continuous function* $\phi : \mathbb{R}^d \to \mathbb{R}$, *we have*

$$(2) \qquad \left|\mathbb{E}_{\tilde{\nu}}\left[\phi\left(\frac{1}{\sqrt{n}}\sum_{k=0}^{n-1} a(X_k)\right)\right] - \mathbb{E}[\phi(N)]\right| \leq \frac{B}{\sqrt{n}} L_{\phi},$$

*where* $N$ *is a* $d$-*dimensional centered Gaussian random variable, centered with covariance matrix* $D(a) = \sum_{k \in \mathbb{Z}}(1-\alpha)^{|k|}\mathbb{E}_{\nu}[a \otimes a \circ T^k]$.

COROLLARY 2.2.6.   *Under Hypothesis* 2.1, *if* $a : \Omega \to \mathbb{R}^d$ *is a* $\nu$-*centered function belonging to* $L^{\infty}(\Omega, \mathbb{R}^d)$, *then, for any Lipschitz continuous function*

$\phi \colon \mathbb{R}^d \to \mathbb{R}$, we have

$$\sup_{s>0, x \in \mathbb{R}^d} \left| \mathbb{E}_{\tilde{\nu}} \left[ \phi \left( x + \varepsilon \sum_{k=0}^{\lfloor s/\varepsilon^2 \rfloor} a(X_k) \right) \right] - \mathbb{E}[\phi(x + \sqrt{s} B_1)] \right|$$

$$\leq \varepsilon L_\phi (B + \|a\|_{L^1(\nu)} + \|B_1\|_{L^1}),$$

where $B_1$ is a $d$-dimensional Gaussian random variable, centered with covariance matrix $D(a)$.

PROOF. Let a real number $s > 0$ be given. If $s < \varepsilon^2$, then we have

$$\left| \mathbb{E}_{\tilde{\nu}} \left[ \phi \left( x + \varepsilon \sum_{k=0}^{\lfloor s/\varepsilon^2 \rfloor} a(X_k) \right) \right] - \phi(x) \right| \leq L_\phi \varepsilon \|a\|_{L^1(\nu)}$$

and

$$|\mathbb{E}[\phi(x + \sqrt{s} B_1)] - \phi(x)| \leq L_\phi \sqrt{s} \|B_1\|_{L^1}.$$

On the other hand, if $s \geq \varepsilon^2$, according to Theorem 2.2.5 [applied to $n = \lfloor \frac{s}{\varepsilon^2} \rfloor + 1$ and to the Lipschitz continuous function $z \mapsto \phi(x + z\varepsilon \sqrt{\lfloor \frac{s}{\varepsilon^2} \rfloor + 1})$], we have

$$\left| \mathbb{E}_{\tilde{\nu}} \left[ \phi \left( x + \varepsilon \sum_{k=0}^{\lfloor s/\varepsilon^2 \rfloor} a(X_k) \right) \right] - \mathbb{E} \left[ \phi \left( x + \varepsilon \sqrt{\left\lfloor \frac{s}{\varepsilon^2} \right\rfloor + 1} B_1 \right) \right] \right| \leq B L_\phi \varepsilon.$$

Moreover, we have

$$\left| \mathbb{E} \left[ \phi \left( x + \varepsilon \sqrt{\left\lfloor \frac{s}{\varepsilon^2} \right\rfloor + 1} B_1 \right) \right] - \mathbb{E}[\phi(x + \sqrt{s} B_1)] \right| \leq L_\phi (\sqrt{s + \varepsilon^2} - \sqrt{s}) \|B_1\|_{L^1}$$

$$\leq L_\phi \frac{\varepsilon^2}{2\sqrt{s}} \|B_1\|_{L^1}$$

$$\leq L_\phi \frac{\varepsilon}{2} \|B_1\|_{L^1}. \qquad \square$$

2.3. *Martingale method.* We recall that the Markov operator $Q_{\alpha,0}$ associated to $(X_n)_{n \in \mathbb{Z}}$ is defined by $Q_{\alpha,0}(f)(\omega) := \mathbb{E}[f(X_{n+1})|X_n = \omega]$, for any $f \in L^1(\Omega, \mathbb{R})$. It is given by the following formula:

$$(3) \qquad Q_{\alpha,0}(f)(\omega) = (1 - \alpha) f \circ T(\omega) + \alpha \mathbb{E}_\nu[f].$$

Using a method introduced by Gordin [17], we get

PROPOSITION 2.3.1. *Let us suppose that Hypothesis* 2.1 *is satisfied. Let* $p \in [1, +\infty]$ *and a function* $f \in L^p(\Omega, \mathbb{C})$ $\nu$-centered. Then, there exists a decomposition of $f$ of the following form:

$$f(X_0) = g_f(X_0, X_{-1}) + h_f(X_0) - h_f(X_{-1}), \qquad \tilde{\nu}\text{-}a.s.,$$

*with* $\mathbb{E}_{\tilde{\nu}}[g_f(X_0, X_{-1})|X_{-1}] = 0$ *and* $h \in L^p(\Omega, \mathbb{C})$.

*Moreover, if* $p \geq 2$, *then we have* $\mathbb{E}_{\tilde{\nu}}[(g_f(X_0, X_{-1}))^2] = D(f)$.

PROOF. Let $p$ and $f$ be as in the hypothesis of the proposition. Let us notice that $Q_{\alpha,0}$ acts continuously on $L^p(\Omega, \mathbb{C})$ and that we have $(Q_{\alpha,0})^n(f) = (1 - \alpha)^n f \circ T^n$, for any integer $n \geq 0$. Using the fact that $X_k \circ \tilde{T} = X_{k+1}$ (for any $k \in \mathbb{Z}$), we can check directly that the functions $g_f(X_0, X_{-1}) = \sum_{n \geq 0}(Q_{\alpha,0})^n(f)(X_0) - \sum_{m \geq 1}(Q_{\alpha,0})^n(f)(X_{-1})$ and $h_f = -\sum_{m \geq 1}(Q_{\alpha,0})^n(f)$ are suitable. Let us prove the second point. By definition, we have

$$D(f) := \lim_{n \to +\infty} \mathbb{E}_{\tilde{\nu}}\left[\left(\frac{1}{\sqrt{n}}\sum_{k=0}^{n-1} f(X_k)\right)^2\right].$$

Moreover, we have

$$\sum_{k=0}^{n-1} f(X_k) = \left(\sum_{k=0}^{n-1} g_f(X_k, X_{k-1})\right) + h_f(X_{n-1}) - h_f(X_{-1}).$$

We conclude by noticing that we have $\mathbb{E}_{\tilde{\nu}}[g_f(X_k, X_{k-1})g_f(X_l, X_{l-1})] = 0$ if $k \neq l$. ☐

The sequence of random variables $(\sum_{k=0}^n g_f(X_k, X_{k-1}))_n$ is a martingale and $h_f(X_0) - h_f(X_{-1}) = h_f \circ X_0 - h_f \circ X_0 \circ \tilde{T}^{-1}$ is a coboundary in $(\tilde{\Omega}, \tilde{\nu}, \tilde{T})$.

This result (Proposition 2.3.1) ensures that, if $a$ is in $L^2(\Omega, \mathbb{R}^d)$, then the sequence of random variables $(a(X_k))_k$ satisfies a central limit theorem. Here, we are interested in a more quantitative question: the rate of convergence in the CLT. To this end, we will use more sophisticated methods (perturbation of operators, Theorem 1.1).

COROLLARY 2.3.2. *Let* $f$ *be a* $\nu$-*centered function belonging to* $L^2(\Omega, \mathbb{R})$ *such that* $D(f) = 0$. *Then we have* $f = 0$, $\nu$-*almost surely.*

PROOF. Let such a function $f$ be given. According to the two previous results, there exists a function $h \in L^2(\Omega, \mathbb{R})$ such that

$$f(X_1) = h(X_1) - h(X_0), \qquad \tilde{\nu}\text{-a.s.}$$

Let us show that $h$ is almost surely constant. Let us suppose that there exists two disjoint measurable subsets $A$ and $B$ of $\Omega$ such that $\nu(h \in A) > 0$ and $\nu(h \in B) > 0$. Then, there exist $\tilde{\omega}_1, \tilde{\omega}_2 \in \tilde{\Omega}$ such that we have

$$X_1(\tilde{\omega}_1) = X_1(\tilde{\omega}_2) = \omega, \qquad h(\omega) \in A, \qquad h(X_0(\tilde{\omega}_1)) \in A, \qquad h(X_0(\tilde{\omega}_2)) \in B$$

and

$$f(X_1(\tilde{\omega}_i)) = h(X_1(\tilde{\omega}_i)) - h(X_0(\tilde{\omega}_i)),$$

and therefore $h(X_0(\tilde{\omega}_1)) = h(\omega) - f(\omega) = h(X_0(\tilde{\omega}_2))$. ☐

2.4. *Proof of Theorem* 2.2.4.   The idea of the method we present here is due to Nagaev [29, 30]. It has been used by many authors (cf., e.g., [19] and [20]). From formula (3), we get, for any integer $n \geq 0$,

$$Q_{\alpha,0}^n(f)(\omega) = (1-\alpha)^n f \circ T^n(\omega) + (1-(1-\alpha)^n)\mathbb{E}_\nu[f].$$

We recall that we have the following relation:

$$Q_{\alpha,0}^n(f)(\omega) = \mathbb{E}_{\tilde{\nu}}[f(X_n)|X_0 = \omega].$$

We will see that the good properties of the Markov operator $Q_{\alpha,0}$ enable us to use the method used in particular in [20].

NOTATION 2.2.   We denote by $\langle \cdot, \cdot \rangle$ the usual scalar product in $\mathbb{R}^d$. If $(\mathcal{B}, \| \cdot \|_\mathcal{B})$ is a complex Banach space, we will use the following notation:

1. We denote by $\mathcal{B}'$ its topological dual (i.e., the set of continuous linear maps from $\mathcal{B}$ in $\mathbb{C}$). We endow this set of the norm $\| \cdot \|_{\mathcal{B}'}$ given by $\|A\|_{\mathcal{B}'} := \sup_{\|f\|_\mathcal{B}=1} |A(f)|$.
2. For any $A \in \mathcal{B}'$ and any $f$ in $\mathcal{B}$, we will use the notation

$$\langle A, f \rangle_* := A(f).$$

3. For any $A \in \mathcal{B}'$, any $g \in \mathcal{B}$, we denote by $g \otimes_* A$ the continuous linear endomorphism of $\mathcal{B}$ defined by

$$(g \otimes_* A)(f) := \langle A, f \rangle_* g.$$

4. We denote by $\mathcal{L}_\mathcal{B}$ the set of continuous linear endomorphisms of $\mathcal{B}$. We endow this set with the norm $\| \cdot \|_{\mathcal{L}_\mathcal{B}}$ given by $\|P\|_{\mathcal{L}_\mathcal{B}} := \sup_{\|f\|_\mathcal{B}=1} \|P(f)\|_\mathcal{B}$.

Let us consider the Banach space $\mathcal{B} := \mathbb{L}^\infty(\Omega, \mathbb{C})$ endowed with the norm $\| \cdot \|_\mathcal{B} = \| \cdot \|_{L^\infty}$. For any $t \in \mathbb{R}^d$, we denote by $Q_{\alpha,t}$ the linear operator on $\mathbb{L}^q(\Omega, \mathbb{C})$ (for any $q \in [1; +\infty]$) defined by

$$Q_{\alpha,t}(f) := Q_{\alpha,0}(e^{i\langle t, a(\cdot)\rangle} f(\cdot)).$$

With this definition, we have

$$(Q_{\alpha,t}(f))(X_n) = \mathbb{E}_{\tilde{\nu}}[e^{i\langle t, a(X_{n+1})\rangle} f(X_{n+1})|X_n].$$

The introduction of these operators is motivated by:

REMARK 2.4.1.   Under Hypothesis 2.1, for any $t \in \mathbb{R}^d$ and any integer $n \geq 1$, we have

$$\mathbb{E}_{\tilde{\nu}}[e^{i\langle t, \sum_{k=0}^{n-1} a(X_k)\rangle}|X_{-1}] = (Q_{\alpha,t})^n(\mathbf{1})(X_{-1}).$$

PROPOSITION 2.4.2. *Let us suppose that Hypothesis 2.1 is satisfied. Let $m \geq 1$ be an integer. If $a \colon \Omega \to \mathbb{R}^d$ is in $L^m(\Omega, \mathbb{R}^d)$, then the function $Q_{\alpha, \cdot} \colon \mathbb{R}^d \to \mathcal{L}_\mathcal{B}$ is $C^m$ on $\mathbb{R}^d$ and, for any integer $k = 1, \ldots, m$ and any $(j_1, \ldots, j_k) \in \{1, \ldots, d\}^k$, we have*

$$\frac{d^k}{dt_{j_1} \cdots dt_{j_k}} Q_{\alpha, t}(f) := i^k Q_{\alpha, t}(a_{j_1}(\cdot) \cdots a_{j_k}(\cdot) f(\cdot))$$

*(where $a_j$ is the $j$th coordinate of $a$).*

We recall that we have $Q_{\alpha, 0}(f) = \alpha \mathbb{E}_\nu[f] + (1 - \alpha) f \circ T$. This can be rewritten as follows:

$$Q_{\alpha, 0}(f) = (\mathbf{1} \otimes_* \nu) f + (1 - \alpha)(f \circ T - \mathbb{E}_\nu[f]).$$

THEOREM 2.4.3 (Perturbation theorem, see [20]). *Let us suppose that Hypothesis 2.1 is satisfied. Let $m \geq 1$ be an integer. We suppose that $a \colon \Omega \to \mathbb{R}^d$ is a $\nu$-centered function belonging to $L^m(\Omega, \mathbb{R}^d)$. Then, there exist a neighborhood $U_0$ of $0$ in $\mathbb{R}^d$ and three nonnegative numbers $c_1, \eta_1, \eta_2$ and four functions $\lambda_{\alpha, \cdot} \in C^m(U_0, \mathbb{C})$, $v_{\alpha, \cdot} \in C^m(U_0, \mathcal{B})$, $\varphi_{\alpha, \cdot} \in C^m(U_0, \mathcal{B}')$ and $N_{\alpha, \cdot} \in C^m(U_0, \mathcal{L}_\mathcal{B})$ such that:*

1. *(Initial values) $\lambda_{\alpha, 0} = 1$, $v_{\alpha, 0} = \mathbf{1}$, $\varphi_{\alpha, 0} = \nu$ and $N_{\alpha, 0}(f) = (1 - \alpha)(f \circ T - \mathbb{E}_\nu[f])$.*
2. *(Initial derivatives) For any $i = 1, \ldots, d$, $\frac{\partial \lambda_{\alpha, t}}{\partial t_1}|_{t=0} = 0$; if $m \geq 2$, then we have $\mathrm{Hess}_t \, \lambda_{\alpha, t}|_{t=0} = -D(a)$, with $D(a) := \sum_{k \in \mathbb{Z}} (1 - \alpha)^{|k|} \mathbb{E}_\nu[a \otimes a \circ T^k]$.*
3. *For any $t$ in $U_0$, we have:*

   (a) *(Decomposition of the operator) For any integer $n \geq 1$, $(Q_{\alpha, t})^n = (\lambda_{\alpha, t})^n v_{\alpha, t} \otimes_* \varphi_{\alpha, t} + (N_{\alpha, t})^n$.*

   (b) *(Dominating eigenvalue) $Q_{\alpha, t} v_{\alpha, t} = \lambda_{\alpha, t} v_{\alpha, t}$, $(Q_{\alpha, t})^* \varphi_{\alpha, t} = \lambda_{\alpha, t} \varphi_{\alpha, t}$ and $\langle \varphi_{\alpha, t}, v_{\alpha, t} \rangle_* = 1$.*

   (c) *$|\lambda_{\alpha, t}| > 1 - \eta_1$.*

   (d) *For any integer $n \geq 1$, we have*

$$\max_{k=0, \ldots, m} \max_{i_1, \ldots, i_k \in \{1, \ldots, d\}} \left\| \frac{d^k}{dt_{i_1} \cdots dt_{i_k}} ((N_{\alpha, t})^n) \right\|_{\mathcal{L}_\mathcal{B}} \leq c_1 (1 - \eta_1 - \eta_2)^n.$$

PROOF. This result is a $d$-dimensional version of Theorem III-8 of [20], page 18. Its proof leads to the implicit function theorem and is exactly the same as the proof of Theorem III-8 of [20]. □

2.4.1. *Regular case.* In our proof, the following theorem plays the same role as the Berry–Esseen lemma in the proof of the rate of convergence in the one-dimensional central limit theorem (see Theorem B of [20], page 12).

PROPOSITION 2.4.4 ([44]). *Let $Q$ be some nondegenerate $d$-dimensional normal distribution. There exist two real numbers $c_0 > 0$ and $\Gamma > 0$ such that, for any real number $T > 0$ and for any Borel probability measure $P$ admitting moments of order $\lfloor \frac{d}{2} \rfloor + 1$, we have*

$$\Pi(P, Q)$$
$$\leq c_0 \left[ \frac{1 + \Gamma}{T} \right.$$
$$\left. + \left( \int_{|t|_\infty < T} \sum_{k=0}^{\lfloor d/2 \rfloor + 1} \sum_{\{i_1, \ldots, i_k\} \in \{1, \ldots, d\}^k} \left| \frac{\partial^k}{\partial t_{i_1} \cdots \partial t_{i_k}} (\varphi_P - \varphi_Q)(t) \right|^2 dt \right)^{1/2} \right].$$

We will prove the following and conclude according to Proposition 2.4.4.

PROPOSITION 2.4.5. *Under the hypotheses of Theorem 2.2.4, if $D(a)$ is invertible, then there exists a real number $\beta > 0$ such that, for any integer $k = 0, \ldots, \lfloor \frac{d}{2} \rfloor + 1$ and any $i_1, \ldots, i_k \in \{1, \ldots, d\}$, we have*

$$\left( \int_{|t|_\infty < \beta \sqrt{n}} \left| \frac{\partial^k}{\partial t_{i_1} \cdots \partial t_{i_k}} (\mathbb{E}_{\tilde{\nu}} [e^{i \langle t, (1/\sqrt{n}) \sum_{l=0}^{n-1} a(X_l) \rangle}] - e^{-(1/2) \langle t, D(a)t \rangle}) \right|^2 dt \right)^{1/2}$$
$$= O\left( \frac{1}{\sqrt{n}} \right).$$

PROOF. The following formula will be useful in the following.

LEMMA 2.4.6. *Let $k$ be a positive integer. Let $b$ be a complex-valued function $C^k$-continuous defined on some open subset $U$ of $\mathbb{R}^d$. Let $n \geq 1$ be an integer. Let us consider the function $u : U \to \mathbb{C}$ given by $u(t) := (b(\frac{t}{\sqrt{n}}))^n$. Then, for any $i_1, \ldots, i_k \in \{1, \ldots, d\}$, we have*

$$\frac{\partial^k}{\partial t_{i_1} \cdots \partial t_{i_k}} u(t)$$
$$= \sum_{\{A_1, \ldots, A_m\} \in Q_k} n(n-1) \cdots (n-m+1) \left( b\left( \frac{t}{\sqrt{n}} \right) \right)^{n-m}$$
$$\times \prod_{i=1}^m \left( \frac{\partial^{\#A_i} b}{\partial t_{i_1^{(i)}} \cdots \partial t_{i_{\#A_i}^{(i)}}} \right) \left( \frac{t}{\sqrt{n}} \right) \frac{1}{n^{k/2}},$$

where $Q_k$ is the set of partitions $\mathcal{A} = \{A_1, \ldots, A_m\}$ of $\{1, \ldots, k\}$ in nonempty subsets $A_i = \{l_1^{(i)}, \ldots, l_{\#A_i}^{(i)}\}$.

Let $c_2 > 0$ and $\beta > 0$ be two real numbers such that the closed ball $\bar{B}_{|\cdot|_\infty}(0, \beta)$ is contained in $U_0$ and such that for any $t \in \bar{B}_{|\cdot|_\infty}(0, \beta)$, we have $|\lambda_{\alpha,t}| \le e^{-c_2\langle t,t\rangle}$ and $e^{-(1/2)\langle t,D(a)t\rangle} \le e^{-c_2\langle t,t\rangle}$. [This is possible because $D(a)$ is invertible and because we have $\mathrm{Hess}_t\,\lambda_{\alpha,t}|_{t=0} = -D(a)$.] In the following, $n$ will be any integer and $t \in \mathbb{R}^d$ any vector satisfying $n \ge 2$ and $|t|_\infty < \beta\sqrt{n}$. For such a couple $(n, t)$, we have $\frac{t}{\sqrt{n}} \in U_0$. Therefore, we have

$$\mathbb{E}_{\tilde{\nu}}[e^{i\langle t,(1/\sqrt{n})\sum_{l=0}^{n-1}a(X_l)\rangle}]$$
$$= \langle \nu, (Q_{\alpha,t/\sqrt{n}})^n \mathbf{1}\rangle_*$$
$$= (\lambda_{\alpha,t/\sqrt{n}})^n \langle \nu, (v_{\alpha,t/\sqrt{n}} \otimes_* \varphi_{\alpha,t/\sqrt{n}})\mathbf{1}\rangle_* + \langle \nu, (N_{\alpha,t/\sqrt{n}})^n \mathbf{1}\rangle_*.$$

1. We start by giving an estimation when $k = 0$. We have

$$\mathbb{E}_{\tilde{\nu}}[e^{i\langle t,(1/\sqrt{n})\sum_{l=0}^{n-1}a(X_l)\rangle}] - e^{-(1/2)\langle t,D(a)t\rangle}$$
$$= (\lambda_{\alpha,t/\sqrt{n}})^n \langle \nu, (v_{\alpha,t/\sqrt{n}} \otimes \varphi_{\alpha,t/\sqrt{n}})\mathbf{1}\rangle_*$$
$$\quad + \langle \nu, (N_{\alpha,t/\sqrt{n}})^n \mathbf{1}\rangle_* - e^{-(1/2)\langle t,D(a)t\rangle}$$
$$= [(\lambda_{\alpha,t/\sqrt{n}})^n - e^{-(1/2)\langle t,D(a)t\rangle}]$$
$$\quad + (\lambda_{\alpha,t/\sqrt{n}})^n (\langle \nu, (v_{\alpha,t/\sqrt{n}} \otimes_* \varphi_{\alpha,t/\sqrt{n}})\mathbf{1}\rangle_* - 1) + \langle \nu, (N_{\alpha,t/\sqrt{n}})^n \mathbf{1}\rangle_*$$
$$= O\left(\frac{1}{\sqrt{n}}|t|_\infty^3 e^{-c_2(1-1/n)\langle t,t\rangle}\right) + O\left(\frac{1}{\sqrt{n}}|t|_\infty e^{-c_2\langle t,t\rangle}\right)$$
$$\quad + c_1(1-\eta_1-\eta_2)^n \frac{|t|_\infty}{\sqrt{n}}$$
$$= O\left(\frac{1}{\sqrt{n}}|t|_\infty^3 e^{-(c_2/2)\langle t,t\rangle}\right) + O\left(\frac{|t|_\infty}{\sqrt{n}}e^{-(c_2/2)\langle t,t\rangle}\right) + c_1(1-\eta_1-\eta_2)^n \frac{|t|_\infty}{\sqrt{n}}.$$

Therefore, we have

$$\left(\int_{|t|_\infty < \beta\sqrt{n}} |\mathbb{E}_{\tilde{\nu}}[e^{i\langle t,(1/\sqrt{n})\sum_{l=0}^{n-1}a(X_l)\rangle}] - e^{-(1/2)\langle t,D(a)t\rangle}|^2\,dt\right)^{1/2} = O\left(\frac{1}{\sqrt{n}}\right).$$

2. Let $k$ be an integer satisfying $1 \le k \le \lfloor\frac{d}{2}\rfloor+1$ and $(i_1, \ldots, i_k) \in \{1, \ldots, d\}^k$. According to Theorem 2.4.3, we have

$$\frac{\partial^k}{\partial t_{i_1}\cdots\partial t_{i_k}}\mathbb{E}_{\tilde{\nu}}[e^{i\langle t,(1/\sqrt{n})\sum_{l=0}^{n-1}a(X_l)\rangle}]$$
$$= \left(\frac{\partial^k}{\partial t_{i_1}\cdots\partial t_{i_k}}((\lambda_{\alpha,t/\sqrt{n}})^n)\right)\langle \nu, (v_{\alpha,t/\sqrt{n}} \otimes_* \varphi_{\alpha,t/\sqrt{n}})\mathbf{1}\rangle_*$$

$$+ O\left( (1 + |t|_\infty^k) \frac{e^{-(c_2/2)\langle t,t\rangle}}{\sqrt{n}} \right) + \frac{1}{n^{k/2}} \left\langle \nu, \frac{\partial^k}{\partial t_{i_1} \cdots \partial t_{i_k}} (N_{\alpha,\cdot})^n \Big|_{t/\sqrt{n}} \mathbf{1} \right\rangle_*$$

$$= \left( \frac{\partial^k}{\partial t_{i_1} \cdots \partial t_{i_k}} ((\lambda_{\alpha,t/\sqrt{n}})^n) \right) \langle \nu, (v_{\alpha,t/\sqrt{n}} \otimes_* \varphi_{\alpha,t/\sqrt{n}}) \mathbf{1} \rangle_*$$

$$\quad + O\left( (1 + |t|_\infty^k) \frac{e^{-(c_2/2)\langle t,t\rangle}}{\sqrt{n}} \right) + \frac{c_1(1 - \eta_1 - \eta_2)^n}{n^{k/2}}$$

$$= \left( \frac{\partial^k}{\partial t_{i_1} \cdots \partial t_{i_k}} ((\lambda_{\alpha,t/\sqrt{n}})^n) \right)$$

$$\quad + O\left( (1 + |t|_\infty^{k+1}) \frac{e^{-(c_2/2)\langle t,t\rangle}}{\sqrt{n}} \right) + \frac{c_1(1 - \eta_1 - \eta_2)^n}{n^{k/2}},$$

since $\langle \nu, (v_{\alpha,t/\sqrt{n}} \otimes_* \varphi_{\alpha,t/\sqrt{n}}) \mathbf{1} \rangle_* - 1 = O(\frac{|t|_\infty}{\sqrt{n}})$ and $\frac{\partial^k}{\partial t_{i_1} \cdots \partial t_{i_k}} ((\lambda_{\alpha,t/\sqrt{n}})^n) = O((1 + |t|_\infty^k)e^{-(c_2/2)\langle t,t\rangle})$. We will estimate the following quantity:

$$\frac{\partial^k}{\partial t_{i_1} \cdots \partial t_{i_k}} ((\lambda_{\alpha,t/\sqrt{n}})^n) - \frac{\partial^k}{\partial t_{i_1} \cdots \partial t_{i_k}} e^{-(1/2)\langle t, D(a)t\rangle}.$$

In the following $b \colon \bar{B}_{|\cdot|_\infty}(0, \beta) \to \mathbb{C}$ will be a function $C^{\lfloor d/2 \rfloor + 1}$ on $\bar{B}_{|\cdot|_\infty}(0, \beta)$ such that $\frac{\partial b}{\partial t_i}(0) = 0$ and $\mathrm{Hess}\, b(0) = -D(a)$ and $|b(t)| \leq e^{-c_2\langle t,t\rangle}$ [we will take $b(t) := \lambda_{\alpha,t}$ and $b(t) := e^{-(1/2)\langle t, D(a)t\rangle}$]. According to Lemma 2.4.6, we have

$$\frac{\partial^k}{\partial t_{i_1} \cdots \partial t_{i_k}} \left( \left( b\left( \frac{t}{\sqrt{n}} \right) \right)^n \right) = \sum_{\mathcal{A} = \{A_1, \ldots, A_m\} \in Q_k} g_{n,m}(\mathcal{A}, b)(t),$$

with

$$g_{n,m}(\mathcal{A}, b)(t) := n(n-1) \cdots (n - m + 1) \left( b\left( \frac{t}{\sqrt{n}} \right) \right)^{n-m}$$

$$\times \prod_{i=1}^m \left( \frac{\partial^{\#A_i} b}{\partial t_{i_1^{(i)}} \cdots \partial t_{i_{\#A_i}^{(i)}}} \right) \left( \frac{t}{\sqrt{n}} \right) n^{-k/2}.$$

For any $\mathcal{A} = \{A_1, \ldots, A_m\} \in Q_k$, we denote by $m_0(\mathcal{A})$ the number of $A_i \in \mathcal{A}$ such that $\#A_i = 1$. Let us notice that we always have $2m \leq m_0(\mathcal{A}) + k$. Indeed, we have

$$k = \sum_{i=1}^m \#A_i \geq m_0(\mathcal{A}) + 2(m - m_0(\mathcal{A})) = 2m - m_0(\mathcal{A}).$$

Therefore, for any $\mathcal{A} = \{A_1, \ldots, A_m\} \in Q_k$, we have

$$|g_{n,m}(\mathcal{A}, b)(t)| \leq n^m e^{-(c_2(n-m)/n)\langle t,t\rangle} O\left( \left( \frac{|t|_\infty}{\sqrt{n}} \right)^{m_0(\mathcal{A})} \right) n^{-k/2}$$

$$= O(n^{(1/2)(2m-(m_0(\mathcal{A})+k))}|t|_\infty^{m_0(\mathcal{A})} e^{-(c_2/2)\langle t,t\rangle})$$

$$= O(|t|_\infty^{m_0(\mathcal{A})} e^{-(c_2/2)\langle t,t\rangle}).$$

(a) If $\mathcal{A} = \{A_1, \ldots, A_m\} \in Q_k$ is such that $2m < m_0(\mathcal{A}) + k$, then, for any $t \in B_{|\cdot|_\infty}(0, \beta\sqrt{n})$, we have

$$|g_{n,m}(\mathcal{A}, b)(t)| = O\left(\frac{|t|_\infty^{m_0(\mathcal{A})}}{\sqrt{n}} e^{-(c_2/2)\langle t,t\rangle}\right).$$

(b) Now, let us consider a partition $\mathcal{A} = \{A_1, \ldots, A_m\} \in Q_k$ such that $2m = m_0(\mathcal{A}) + k$. Then, $\mathcal{A}$ is made of subsets of $\{1, \ldots, d\}$ containing at most two elements. For such a partition $\mathcal{A}$, for any $t \in B_{|\cdot|_\infty}(0, \beta\sqrt{n})$, we have

$$|g_{n,m}(\mathcal{A}, \lambda_{\alpha,\cdot})(t) - g_{n,m}(\mathcal{A}, e^{-(1/2)\langle \cdot, D(a)\cdot\rangle})(t)|$$

$$= O\left(\frac{1}{\sqrt{n}}(1 + |t|_\infty^{m_0(A)+3}) e^{-(c_2/2)\langle t,t\rangle}\right).$$

Indeed, we have

$$\frac{\partial}{\partial t_i}(\lambda_{\alpha,\cdot} - e^{-(1/2)\langle \cdot, D(a)\cdot\rangle})\left(\frac{t}{\sqrt{n}}\right) = O\left(\frac{|t|_\infty^2}{n}\right),$$

$$(\lambda_{\alpha,t/\sqrt{n}})^{n-m} - e^{-((n-m)/2n)\langle t,D(a)t\rangle} = O\left(\frac{|t|_\infty^3}{\sqrt{n}} e^{-(c_2/2)\langle t,t\rangle}\right)$$

and

$$\frac{\partial^2}{\partial t_i \partial t_j}(\lambda_{\alpha,\cdot} - e^{-(1/2)\langle \cdot, D(a)\cdot\rangle})\left(\frac{t}{\sqrt{n}}\right) = O\left(\frac{|t|_\infty}{\sqrt{n}}\right). \qquad \square$$

2.4.2. *Degenerate case.* In this section we suppose that the matrix $D(a)$ is degenerate (i.e., noninvertible). There exists a matrix $A \in GL(\mathbb{R}^d)$ such that we have

$$A \cdot D(a) \cdot {}^T A = J_l := \begin{pmatrix} I_l & 0_{l,d-l} \\ 0_{d-l,l} & 0_{d-l,d-l} \end{pmatrix},$$

where $l$ is the rank of $D(a)$, $I_l$ is the $l$-dimensional identity matrix and $0_{m,n}$ is the $(m,n)$-dimensional null matrix. By replacing function $a(\cdot)$ by $A \cdot a(\cdot)$, we can (and we will) assume that we have $D(a) = J_l$. According to the previous subsection, we have

$$(4) \quad \Pi\left(\tilde{\nu}_*\left(\left(\frac{1}{\sqrt{n}}\sum_{k=0}^{n-1} a_1(X_k), \ldots, \frac{1}{\sqrt{n}}\sum_{k=0}^{n-1} a_l(X_k)\right)\right), \mathcal{N}(0, I_l)\right) = O(\varepsilon),$$

where $a_i$ is the $i$th coordinate of $a$. For any $i = l+1, \ldots, d$, we have $D(a_i) = 0$ and therefore, according to Corollary 2.3.2, we have $a_i = 0$ almost surely.

2.5. *Proof of Theorem* 2.2.5. The sequence of random variables $(a(X_k))_k$ satisfies the hypotheses of Theorem 1.1. Indeed, for any integers $\alpha, \beta, \gamma \geq 0$ satisfying $1 \leq \alpha + \beta + \gamma \leq 3$, any integers $1 \leq k \leq k+p \leq k+p+q \leq k+p+l$, and any integers $i_1, i_2, i_3$, we have

$$\mathbb{E}_{\tilde{\nu}}[(a_{i_1}(X_{k+p}))^\alpha (a_{i_2}(X_{k+p+q}))^\beta (a_{i_3}(X_{k+p+l}))^\gamma | X_0, X_1, \ldots, X_k]$$
$$= Q_{\alpha,0}^p(\psi)(X_k)$$

and

$$\mathbb{E}_{\tilde{\nu}}[(a_{i_1}(X_{k+p}))^\alpha (a_{i_2}(X_{k+p+q}))^\beta (a_{i_3}(X_{k+p+l}))^\gamma] = \mathbb{E}_\nu[\psi],$$

with $\psi = \psi_{q,l,\alpha,\beta,\gamma,i_1,i_2,i_3} := a_{i_1}^\alpha \times Q_{\alpha,0}^q(a_{i_2}^\beta \times (Q_{\alpha,0}^{l-q} a_{i_3}^\gamma))$. We get

$$|\mathrm{Cov}(G(X_0, \ldots, X_k), (a_{i_1}(X_{k+p}))^\alpha (a_{i_2}(X_{k+p+q}))^\beta (a_{i_3}(X_{k+p+l}))^\gamma)|$$
$$\leq \|Q_{\alpha,0}^p(\psi - \mathbb{E}_\nu[\psi])\|_{L^\infty} \|G(X_0, \ldots, X_k)\|_{L^1}$$
$$\leq (1-\alpha)^p \|\psi - \mathbb{E}_\nu[\psi]\|_{L^\infty} \|G(X_0, \ldots, X_k)\|_{L^1}$$
$$\leq 2(1-\alpha)^p \|\psi\|_{L^\infty} \|G(X_0, \ldots, X_k)\|_{L^1}$$
$$\leq 2(1-\alpha)^p \|a\|_\infty^{\alpha+\beta+\gamma} \|G(X_0, \ldots, X_k)\|_{L^1}$$
$$\leq 2(1-\alpha)^p (1 + \|a\|_\infty^3) \|G(X_0, \ldots, X_k)\|_{L^1}.$$

## 3. Application to the Sinai billiard.

The billiard system considered here has been studied in many articles since the fundamental article of Sinai [40]. Let us mention [8, 9, 10, 11, 24]. The question of the CLT in this context has been studied in many articles [9, 10, 12, 43].

Here, we are interested in the question of speed of convergence in the CLT.

A first result has been established in [33] for one-dimensional observables (the speed is estimated in the sense of the uniform norm of the difference between repartition functions). This result has been extended in [34] (for $d$-dimensional observables, the speed being estimated in the sense of the Prokhorov metric).

The speed obtained in these two papers is in $n^{-(1/2)+\alpha}$ for all $\alpha > 0$. Here, we establish a rate of convergence in the CLT in $n^{-1/2}$ in the sense of the Kantorovich metric. This result is an application of Theorem 1.1.

3.1. *The model.* We are interested in the behavior of a point particle moving with unit speed in some domain $Q$ of the torus $\mathbb{T}^2$, the complement of which is a finite union of open sets $\mathcal{O}_1, \ldots, \mathcal{O}_I$ called obstacles. Each obstacle $\mathcal{O}_i$ is a strictly convex open set, the boundary of which is $C^3$ and the curvature of the boundary is never null. See Figure 2.

We suppose that the closures of the obstacles are pairwise disjoint. We suppose that the point particle moves in $Q$ with unit speed and elastic reflection off the obstacles. See Figure 3.

**3.2.** *The billiard flow.* Let us notice that, when the particle hits an obstacle, the couple position-speed is ambiguously defined: incoming and outgoing vectors coexist. To avoid this problem, we decide to take the following convention: when a particle hits an obstacle, its position-speed couple corresponds to the outgoing vector. Let us be more precise. For all $q$ in $\partial Q$, we denote by $\vec{n}(q)$ the unit vector normal to $\partial Q$ in $q$ directed to the interior of $Q$.

The set of configurations is the set $Q_1$ given by

$$Q_1 := \{(q, \vec{v}) : q \in Q, \vec{v} \in T_q Q, \|\vec{v}\| = 1 \text{ and } (q \in \partial Q \Rightarrow \langle \vec{n}(q), \vec{v} \rangle \geq 0)\}.$$

We call billiard flow the flow $(Y_t)_t$ defined on $Q_1$ such that, for all $t > 0$, all $(q, \vec{v}) \in Q_1$ and all $(q', \vec{v}') \in Q_1$, the fact that $Y_t(q, \vec{v}) = (q', \vec{v}')$ means that "if a particle is at $q$ with the speed $\vec{v}$ at time 0, then it will be at $q'$ with speed $\vec{v}'$ at time $t$."

This flow preserves the normalized Lebesgue $\mu$ on $Q_1$.

According to the description of our model, it is natural to study the model corresponding to the times when the particle hits an obstacle [cf. the system $(M, \nu, T)$ below].

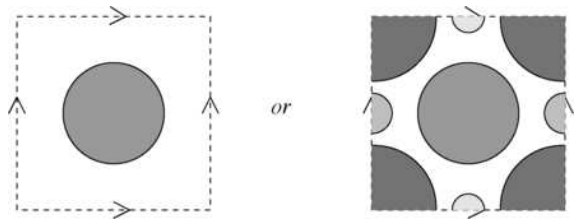

Fig. 2.

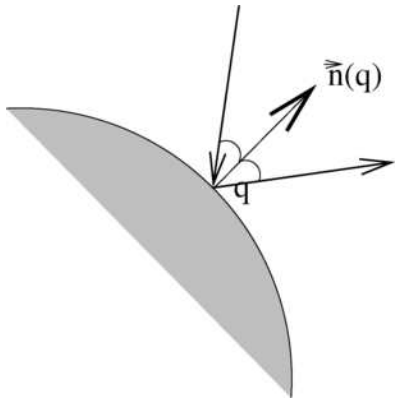

Fig. 3.

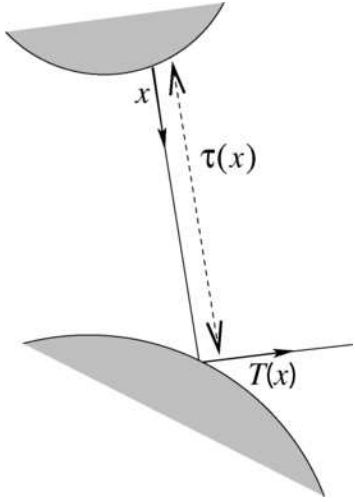

Fig. 4.

3.3. *The billiard transformation.* Let us consider the billiard system $(M, \nu, T)$ defined as follows:

(a) $M$ is the set of configurations of $Q_1$ corresponding to the times when the particle meets an obstacle, that is,

$$M := \{(q, \vec{v}) : q \in \partial Q, \vec{v} \in T_q Q, \|\vec{v}\| = 1, \langle \vec{n}(q), \vec{v} \rangle \geq 0\}.$$

(b) For any $i = 1, \ldots, I$, we write $l_i$ the length of the boundary $\partial \mathcal{O}_i$ of the obstacle $\mathcal{O}_i$. We parametrize $M$ by $G : M \to \bigcup_{i=1}^{I} (\{i\} \times \frac{\mathbb{R}}{l_i \mathbb{Z}} \times [-\pi/2; \pi/2])$ defined by $G(q, \vec{v}) = (i, r, \varphi)$ if $q \in \mathcal{O}_i$, if $r$ is the curvilinear abscissa of $q$ on $\mathcal{O}_i$, and $\varphi$ is the angular measure taken in $[-\pi/2; \pi/2]$ of the angle between $\vec{n}(q)$ and $\vec{v}$.

(c) $\nu$ is the Borel probability measure on $M$ of the following form:

$$\nu(A) = \frac{1}{C} \sum_{i=1}^{I} \int_{\{(r, \varphi) : G^{-1}(i, r, \varphi) \in A\}} \cos(\varphi) \, dr \, d\varphi,$$

where $C$ is some constant.

(d) $T$ is the transformation of $M$ that, at the configuration $(q, \vec{v}) \in M$ of a particle at the time just after a reflection, associates the configuration $(q', \vec{v}') \in M$ at the time just after the following reflection off $\partial Q$ (cf. Figure 4).

(e) We also define the function $\tau : M \to [0; +\infty[$ where $\tau(q, \vec{v})$ is the time to wait for a particle at $q$ with speed $\vec{v}$ until the next reflection off $\partial Q$ (cf. Figure 4):

$$\tau(q, \vec{v}) := \min\{t > 0 : q + t\vec{v} \in \partial Q\}.$$

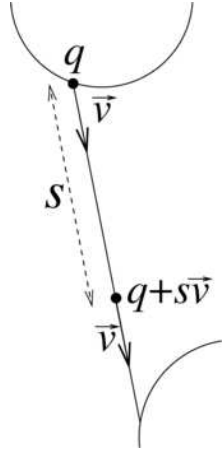

Fig. 5.

Let us specify the link between billiard transformation and billiard flow. The flow $(Y_t)_t$ can be viewed as the special flow over the dynamical system $(M, \nu, T)$ associated to the roof-function $\tau$. This is very natural: we identify $((q, \vec{v}), s)$ with $(q + s\vec{v}, \vec{v})$ (cf. Figure 5).

In the following, we suppose that the billiard system has *finite horizon*, that is, that function $\tau$ is uniformly bounded.

In Figure 2 only the second domain corresponds to a billiard system with finite horizon.

3.4. *About the regularity of $T$.* Let $R_0$ be the set of configurations corresponding to a vector tangent to an obstacle:

$$R_0 := \{(q, \vec{v}) \in M : \langle \vec{n}(q), \vec{v} \rangle = 0\}.$$

The study of the billiard is complicated by the discontinuity of $T$ at points of $T^{-1}(R_0)$ (cf. Figure 6).

However, we know that, for any integer $k \geq 1$, the transformation $T^k$ defines a $C^1$-diffeomorphism from $M \setminus \bigcup_{j=0}^{k} T^{-j}(R_0)$ onto $M \setminus \bigcup_{j=0}^{k} T^{j}(R_0)$. Moreover, the sets $\bigcup_{j=0}^{k} T^{-j}(R_0)$ and $\bigcup_{j=0}^{k} T^{j}(R_0)$ are finite union of $C^1$-curves.

3.5. *Hyperbolic properties of the billiard transformation.* For any $C^1$-curve $\gamma$ of $M$, we define

$$l(\gamma) := \int_\gamma \sqrt{dr^2 + d\varphi^2},$$

using the parametrization of $M$ by the function $G$ previously defined.

PROPOSITION 3.5.1.   *There exist two real numbers $C_0 > 0$ and $\alpha_0 \in \,]0;1]$ such that, for all $x$, there exist $C^1$-curves $\gamma^s(x)$ (stable curve) and $\gamma^u(x)$ (unstable curve) of $M$ containing $x$ (with positive length for $\nu$-almost every $x$) such that, for any integer $n \geq 0$, all $y, z \in \gamma^s(x)$ and all $y', z' \in \gamma^u(x)$, we have*

$$d(T^n(y), T^n(z)) \leq C_0 \alpha_0^n \sqrt{d(y,z)}$$

*and*

$$d(T^{-n}(y'), T^{-n}(z')) \leq C_0 \alpha_0^n \sqrt{d(y', z')}.$$

3.6. *The functional sets $\mathcal{H}_{\eta,m}$.*   Because of the discontinuities of $T$, if $\phi \colon M \to \mathbb{R}$ is a Hölder continuous function, then $\phi \circ T^m$ is generally not a Hölder continuous function. This observation leads us to the introduction of the sets $\mathcal{H}_{\eta,m}$ defined below. These spaces will be such that, if $f$ is $\eta$-Hölder continuous, then $f \circ T^m$ is in $\mathcal{H}_{\eta,m}$.

Let a real number $\eta \in \,]0;1]$ be given. For any $m$, we consider the set $\mathcal{H}_{\eta,m}$ of bounded functions $\phi \colon M \to \mathbf{C}$ such that the following quantity is finite:

$$C_\phi^{(\eta,m)} := \sup_{C \in \mathcal{C}_m} \sup_{x,y \in C, x \neq y} \frac{|\phi(x) - \phi(y)|}{(\max(d(x,y), \ldots, d(T^m(x), T^m(y))))^\eta},$$

where $\mathcal{C}_m$ is the set of the connected components of $M \setminus \bigcup_{j=0}^m T^{-j}(R_0)$ and $d$ is the metric defined on each connected component of $M$ by $d((q, \vec{v}), (q', \vec{v'})) = \sqrt{|r - r'| + |\varphi - \varphi'|^2}$ if $G(q, \vec{v}) = (i, r, \varphi)$ and $G(q', \vec{v'}) = (i, r', \varphi')$.

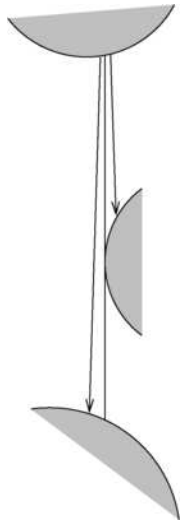

Fɪɢ. 6.

The set $\mathcal{H}_{\eta,m}$ can be understood as the set of functions that are Hölder continuous in the $m$ future configurations. These classes of functions have been introduced in [33].

Let us notice that the function $\tau$ is in $\mathcal{H}_{1,1}$.

In the following section, we give a decorrelation result for these classes of functions. Before recalling this result, let us make some comments about the classes of functions $\mathcal{H}_{\eta,m}$.

PROPOSITION 3.6.1. *Let a real number* $\eta \in ]0;1]$ *and an integer* $m_0 \geq 1$ *be given. For any functions* $\phi$ *and* $\psi$ *belonging to* $\mathcal{H}_{\eta,m_0}$, *we have:*

1. *The functions* $\phi$ *and* $\psi$ *are uniformly bounded [because the set* $M \setminus \bigcup_{j=0}^{m_0} T^{-j}(R_0)$ *has only a finite number of connected components].*
2. *The function* $\phi + \psi$ *is in* $\mathcal{H}_{\eta,m_0}$ *and we have*

$$C_{\phi+\psi}^{(\eta,m_0)} \leq C_{\phi}^{(\eta,m_0)} + C_{\psi}^{(\eta,m_0)}.$$

3. *The product* $\phi \cdot \psi$ *is in* $\mathcal{H}_{\eta,m_0}$ *and we have*

$$C_{\phi\cdot\psi}^{(\eta,m_0)} \leq C_{\phi}^{(\eta,m_0)} \|\psi\|_\infty + C_{\psi}^{(\eta,m_0)} \|\phi\|_\infty.$$

4. *For any integer* $m \geq 0$, $\phi \circ T^m$ *is in* $\mathcal{H}_{\eta,m_0+m}$ *and we have*

$$C_{\phi\circ T^m}^{(\eta,m_0+m)} \leq C_{\phi}^{(\eta,m_0)}.$$

5. *For any integer* $m \geq 0$, *the function* $\phi$ *is in* $\mathcal{H}_{\eta,m_0+m}$ *and we have*

$$C_{\phi}^{(\eta,m_0+m)} \leq C_{\phi}^{(\eta,m_0)}.$$

PROOF OF POINT 4. Let $x$ and $y$ be two points of $M$ belonging to the same connected component of $M \setminus \bigcup_{j=0}^{m+m_0} T^{-j}(R_0)$. Then $T^m(x)$ and $T^m(y)$ belong to the same connected component of $M \setminus \bigcup_{j=0}^{m_0} T^{-j}(R_0)$ and we have

$$|\phi(T^m(x)) - \phi(T^m(x))|$$
$$\leq C_{\phi}^{(\eta,m_0)} \max(d(T^m(x), T^m(y)), \dots, d(T^{m+m_0}(x), T^{m+m_0}(y)))^\eta. \quad \square$$

3.7. *A decorrelation property.* The following result has been established in [33] (cf. Proposition 1.2 and Corollary B.2 of [33]) with the use of the method developed by Young in [43].

PROPOSITION 3.7.1. *Let a real number* $\eta \in ]0;1]$ *be given. For any real number* $R > 1$, *there exist two real numbers* $C_{\eta,R} > 0$ *and* $\delta_{\eta,R} \in ]0;1[$ *such that, for any integers* $m_1 \geq 0$ *and* $m_2 \geq 0$, *for all* $\phi \in \mathcal{H}_{\eta,m_1}$ *and* $\psi \in \mathcal{H}_{\eta,m_2}$, *for any integer* $n \geq 0$, *we have*

(5)
$$|\operatorname{Cov}(\phi, \psi \circ T^n)|$$
$$\leq C_{\eta,R}(\|\phi\|_\infty \|\psi\|_\infty + C_{\phi}^{(\eta,m_1)} \|\psi\|_\infty + \|\phi\|_\infty C_{\psi}^{(\eta,m_2)}) \delta_{\eta,R}^{n-Rm_1}.$$

We will not use directly this proposition; we will use a slight modification of it: reading the proof of Theorem B.1 of [33], we can notice that, in formula (5), coefficient $C_\psi^{(\eta, m_2)}$ can be replaced by the regularity coefficient of $\psi$ on the stable curves:

$$C_\psi^{(\eta, (s))} := \sup_{x \in M} \sup_{y, z \in \gamma^s(x)} \frac{|\psi(y) - \psi(z)|}{d(y, z)^\eta}$$

[by replacing, in the proof of Theorem B.1, the definition of $\hat{\psi}_k(x)$ by the infimum of $\psi \circ \tilde{T}_d^k$ on the stable curve containing $x$].

### 3.8. Theorem.

THEOREM 3.8.1.   Let a real number $\eta \in \, ]0; 1]$ and an integer $m_0 \geq 1$ be given. Let $f : M \to \mathbb{R}^d$ be a bounded function, the coordinates of which are in $\mathcal{H}_{\eta, m_0}$. For any $k$, we write $Y_k := f \circ T^k$ and $S_k := Y_1 + \cdots + Y_k$. Then, the following limit exists:

$$\Sigma^2 := \lim_{n \to +\infty} \frac{1}{n} (\mathbb{E}[S_n^{\otimes 2}]).$$

If $\Sigma^2 = 0$, then $(S_n)_n$ is bounded in $L^2$.

Otherwise, the sequence of random variables $(\frac{S_n}{\sqrt{n}})_{n \geq 1}$ converges in distribution to a Gaussian random variable $N$ with null expectation and covariance matrix $\Sigma^2$ and there exists a real number $B > 0$ such that, for all $n \geq 1$ and all Lipschitz continuous function $\phi : \mathbb{R}^d \to \mathbb{R}$, we have

$$\left| \mathbb{E}\left[ \phi\left( \frac{S_n}{\sqrt{n}} \right) \right] - \mathbb{E}[\phi(N)] \right| \leq \frac{B L_\phi}{\sqrt{n}}.$$

PROOF.   For all $n$, $Y_1 + \cdots + Y_n$ has the same distribution as $X_1 + \cdots + X_n$ with $X_k := f \circ T^{-k}$. Therefore, it suffices to show that the sequence $(X_n)_n$ satisfies the hypotheses of Theorem 1.1. Let us take $M := \max(1, \|f\|_\infty)$. Let $a, b, c$ be integers satisfying $a, b, c \geq 1$ and $1 \leq a + b + c \leq 3$. Let $i, j, k, p, q, l$ be integers satisfying $1 \leq i \leq j \leq k \leq k + p \leq k + p + q \leq k + p + l$. Let $i_1, i_2, i_3$ be integers belonging to $\{1, \ldots, d\}$. Let $F : \mathbb{R}^d \times ([-M; M]^d)^3 \to \mathbb{R}$ be a bounded, differentiable function, with bounded differential. We have

$$|\text{Cov}(F(S_{i-1}, X_i, X_j, X_k), (X_{k+p}^{(i_1)})^a (X_{k+p+q}^{(i_2)})^b (X_{k+p+l}^{(i_3)})^c)|$$
$$= |\text{Cov}(F(f \circ T^{-1} + \cdots + f \circ T^{-(i-1)}, f \circ T^{-i}, f \circ T^{-j}, f \circ T^{-k}),$$
$$f_{i_1}^a \circ T^{-(k+p)} f_{i_2}^b \circ T^{-(k+p+q)} f_{i_3}^c \circ T^{-(k+p+l)})|$$
$$= |\text{Cov}(\phi_0, \psi_0 \circ T^{p+l})|,$$

where we define

$$\phi_0 := f_{i_3}^c f_{i_2}^b \circ T^{l-q} f_{i_1}^a \circ T^l$$

and

$$\psi_0 := F(f \circ T^{k-i+1} + \cdots + f \circ T^{k-1}, f \circ T^{k-i}, f \circ T^{k-j}, f).$$

Let us use formula (5) modified with $C_\psi^{(\eta,(s))}$ instead of $C_\psi^{(\eta,m_2)}$.

LEMMA 3.8.2. *The function $\psi_0$ is in $\mathcal{H}_{\eta/2,m_0+k-i+1}$ and we have*

$$C_{\psi_0}^{(\eta/2,(s))} \le \|DF\|_\infty C_f^{(\eta,m_0)} \frac{C_0^\eta}{1-\alpha_0^\eta}.$$

PROOF.   Let three points $x, y, z$ in $M$ be such that $y$ and $z$ are in $\gamma^s(x)$. Then we have

$$|f(y) - f(z)| \le C_f^{(\eta,m_0)} C_0^\eta d(y,z)^{\eta/2},$$

$$|f(T^{k-i}(y)) - f(T^{k-i}(z))| \le C_f^{(\eta,m_0)}(C_0\alpha_0^{k-i})^\eta d(y,z)^{\eta/2},$$

$$|f(T^{k-j}(y)) - f(T^{k-j}(z))| \le C_f^{(\eta,m_0)}(C_0\alpha_0^{k-j})^\eta d(y,z)^{\eta/2}$$

and

$$|f \circ T^{k-i+1}(y) + \cdots + f \circ T^{k-1}(y) - (f \circ T^{k-i+1}(z) + \cdots + f \circ T^{k-1}(z))|$$

$$\le \sum_{m=1}^{i-1} |f \circ T^{k-i+m}(y) - f \circ T^{k-i+m}(z)|$$

$$\le \sum_{m=1}^{i-1} C_f^{(\eta,m_0)}(C_0\alpha_0^{k-i+m})^\eta d(y,z)^{\eta/2}$$

$$\le C_f^{(\eta,m_0)} C_0^\eta \frac{\alpha_0^{(k-i+1)\eta}}{1-\alpha_0^\eta} d(y,z)^{\eta/2}. \qquad \square$$

Since $f$ is in $\mathcal{H}_{\eta,m_0}$, the function $\phi_0$ is in $\mathcal{H}_{\eta/2,m_0+l}$ and we have

$$\|\phi_0\|_\infty \le \|f\|_\infty^{a+b+c} \le (1 + \|f\|_\infty^3)$$

and

$$C_{\phi_0}^{(\eta/2,m_0+l)} \le (a+b+c)C_f^{(\eta/2,m_0)}\|f\|_\infty^{a+b+c-1} \le 3C_f^{(\eta/2,m_0)}(1 + \|f\|_\infty^2).$$

Therefore, according to (5) modified with $C_\psi^{(\eta,(s))}$ instead of $C_\psi^{(\eta,m_2)}$, we have

$$|\text{Cov}(\phi_0, \psi_0 \circ T^{p+l})|$$

$$\le C_{\eta/2,R}(\|\phi_0\|_\infty + C_{\phi_0}^{(\eta/2,m_0+l)})(\|\psi_0\|_\infty + C_{\psi_0}^{(\eta/2,(s))})\delta_{\eta/2,R}^{(p+l)-R(m_0+l)}$$

$$\le K(\|F\|_\infty + \|DF\|_\infty)\varphi_{p,l},$$

with

$$K := C_{\eta/2,R}(1 + 3C_f^{(\eta/2,m_0)})(1 + \|f\|_\infty^3)\left(1 + C_f^{(\eta,m_0)}\frac{C_0^\eta}{1 - \alpha_0^\eta}\right)\delta_{\eta/2,R}^{-Rm_0}$$

and

$$\varphi_{p,l} := \delta_{\eta/2,R}^{p - (R-1)l}.$$

We have

$$\sum_{p \geq 1} p \max_{l=0,\ldots,\lfloor p/\lfloor R \rfloor \rfloor} \varphi_{p,l} \leq \sum_{p \geq 1} p\delta_{\eta/2,R}^{p - ((R-1)/\lfloor R \rfloor)p}$$

$$= \sum_{p \geq 1} p\delta_{\eta/2,R}^{((\lfloor R \rfloor + 1 - R)/\lfloor R \rfloor)p} < +\infty. \qquad \square$$

## APPENDIX

PROOF OF THEOREM 1.1. Our proof of Theorem 1.1 follows the same scheme as [25].

First, since $(X_k)_{k \geq 1}$ is stationary, we notice that, for any integer $n \geq 1$, we have

$$(6) \quad \mathbb{E}\left[\left(\frac{S_n}{\sqrt{n}}\right)^{\otimes 2}\right] = \mathbb{E}[X_1^{\otimes 2}] + \sum_{k=1}^{n-1}\left(1 - \frac{k}{n}\right)(\mathbb{E}[X_1 \otimes X_{k+1}] + \mathbb{E}[X_{k+1} \otimes X_1]).$$

Therefore, according to (1), $\Sigma^2$ exists and we have

$$\Sigma^2 = \mathbb{E}[X_1^{\otimes 2}] + \sum_{k \geq 1}(\mathbb{E}[X_1 \otimes X_{k+1}] + \mathbb{E}[X_{k+1} \otimes X_1]).$$

Moreover, according to (6), we have

$$|\mathbb{E}[S_n^{\otimes 2}] - n\Sigma^2|_\infty \leq 2\sum_{k \geq 1} k|\mathbb{E}[X_1 \otimes X_{k+1}]|_\infty \leq 4CM\sum_{k \geq 1} k\varphi_{k,0}.$$

Hence, if $\Sigma^2 = 0$, then the sequence of random variables $(S_n)_{n \geq 1}$ is bounded in $L^2(\Omega, \mathbb{R}^d)$.

Let us now suppose that $\Sigma^2$ is nonnull. Then, there exists $k \in \{1, \ldots, d\}$ and a $d$-dimensional orthogonal matrix $O$ such that $O \cdot \Sigma^2 \cdot O^{-1}$ is diagonal. Therefore, there exists an invertible matrix $A$ such that

$$A \cdot \Sigma^2 \cdot {}^tA = J_k,$$

where ${}^tA$ denotes the matrix transposed to $A$ and where $J_k$ is the $d$-dimensional diagonal matrix such that the first $k$ diagonal elements are equal to 1 and the others to 0.

In the following, we will suppose that $\Sigma^2$ is the $d$-dimensional identity matrix $I_d$. This is not a restrictive hypothesis: it suffices to replace $d$ by $k$ and $(X_n)_n$ by $(\tilde{X}_n := (\tilde{X}_n^{(1)}, \ldots, \tilde{X}_n^{(k)}))_n$ where $\tilde{X}_n^{(i)}$ is the $i$th coordinate of $\tilde{X}_n := A \cdot X_n$ [since the random variables $(\tilde{X}_1^{(j)} + \cdots + \tilde{X}_n^{(j)})_{n \geq 1}$ are bounded in $L^2(\Omega, \mathbb{R})$ for all $j = k+1, \ldots, d$ and the norms on $\mathbb{R}^d$ are equivalent].

As in [37], we will use an inductive proof. The idea is to prove the existence of a real number $A \geq 1$ such that the following property $(\mathcal{P}_n(A))$ is satisfied for any integer $n \geq 2$:

$$(\mathcal{P}_n(A)) : \forall k = 1, \ldots, n-1, \forall \phi \in \mathbf{Lip}(\mathbb{R}^d, \mathbb{R})$$

$$|\mathbb{E}[\phi(S_k)] - \mathbb{E}[\phi(\sqrt{k}N)]| \leq AL_\phi,$$

where $N$ is a $d$-dimensional Gaussian random variable with expectation $0$ and covariance matrix $I_d$. Let us define $V_n := \mathbb{E}[S_n^{\otimes 2}]$ and $v_n := \mathbb{E}[S_n^{\otimes 2}] - \mathbb{E}[S_{n-1}^{\otimes 2}]$. We have

$$v_n = \mathbb{E}[X_1^{\otimes 2}] + \sum_{k=1,\ldots,n-1} (\mathbb{E}[X_1 \otimes X_{k+1}] + \mathbb{E}[X_{k+1} \otimes X_1]).$$

Hence $(v_n)_{n \geq 1}$ converges to $\Sigma^2 = I_d$. There exists $n_0 \geq 1$ such that for any integer $n \geq n_0$, the eigenvalues of $v_n$ are between $\frac{1}{2}$ and $\frac{3}{2}$. In the following, we will suppose the existence of a sequence $(N_i)_{i \geq 0}$ of independent identically distributed Gaussian random variables with expectation $0$ and covariance matrix $I_d$ such that $(N_i)_{i \geq 0}$ is independent of $(X_k)_{k \geq 0}$. The main part of the proof is to establish the following result.

PROPOSITION A.1. *Under the hypotheses of Theorem* 1.1, *there exist a real number $K \geq M$ and a continuous decreasing function $\psi : [1, +\infty) \to ]0; +\infty[$ satisfying $\lim_{\varepsilon \to +\infty} \psi(\varepsilon) = 0$ such that for any integer $n \geq 9n_0$ and any real number $A \geq M$, if we have $(\mathcal{P}_n(A))$, then, for any real number $\varepsilon \geq 1$ and any Lipschitz continuous function $\phi : \mathbb{R}^d \to \mathbb{R}$, we have*

$$|\mathbb{E}[\phi(S_n + \varepsilon Y)] - \mathbb{E}[\phi(S_{n_0-1} + T_{n_0-1,n} + \varepsilon Y)]| \leq K(1 + A\psi(\varepsilon))L_\phi,$$

*where $Y$ and $T_{n_0-1,n}$ are two $\nu$-centered Gaussian random variables independent of $(X_k)_{k \geq 0}$, with covariance matrices $I_d$ and $V_n - V_{n_0-1}$, respectively.*

Let us show how we can conclude once this result is proved.

1. Let us write $A_1 := d\sqrt{9n_0} + \max_{m=0,\ldots,9n_0} \|S_m\|_{L_1}$. For any $A \geq A_1$, property $(\mathcal{P}_{9n_0}(A))$ is satisfied.

2. Let us show that there exists a real number $A_0 \geq M$ such that, for any integer $n \geq 9n_0$ and any real number $A \geq A_0$, we have $(\mathcal{P}_n(A)) \Rightarrow (\mathcal{P}_{n+1}(A))$.

Let an integer $n \geq 9n_0$ and a real number $A \geq M$ be given such that property $(\mathcal{P}_n(A))$ is satisfied. Let $\phi \colon \mathbb{R}^d \to \mathbb{R}$ be any Lipschitz continuous function. Then, according to Proposition A.1, we have

$$|\mathbb{E}[\phi(S_n + \varepsilon Y)] - \mathbb{E}[\phi(S_{n_0-1} + T_{n_0-1,n} + \varepsilon Y)]| \leq K(1 + A\psi(\varepsilon))L_\phi.$$

Since we have

$$|\mathbb{E}[\phi(S_n + \varepsilon Y)] - \mathbb{E}[\phi(S_n)]| \leq L_\phi \varepsilon \mathbb{E}[|Y|_\infty]$$

and

$$|\mathbb{E}[\phi(S_{n_0-1} + T_{n_0-1,n} + \varepsilon Y)] - \mathbb{E}[\phi(T_{n_0-1,n})]|$$
$$\leq L_\phi(\varepsilon \mathbb{E}[|Y|_\infty] + \mathbb{E}[|S_{n_0-1}|_\infty]),$$

we get

$$|\mathbb{E}[\phi(S_n)] - \mathbb{E}[\phi(T_{n_0-1,n})]| \leq K_0(1 + A\psi(\varepsilon) + 2\varepsilon)L_\phi,$$

with $K_0 := K + \mathbb{E}[|Y|_\infty] + \mathbb{E}[|S_{n_0-1}|_\infty]$. Let us now estimate the following quantity:

$$|\mathbb{E}[\phi(\sqrt{n}Y)] - \mathbb{E}[\phi(T_{n_0-1,n})]|.$$

Let $O = O_{n,n_0}$ be an orthogonal matrix such that $O(V_n - V_{n_0-1})O^{-1}$ is diagonal with nonnegative diagonal coefficients. Let us denote by $\Delta_{n,n_0}$ the diagonal matrix with nonnegative diagonal coefficients such that

$$(\Delta_{n,n_0})^2 = O(V_n - V_{n_0-1})O^{-1}.$$

Let us define $M_{n,n_0}$ as follows:

$$M_{n,n_0} := O^{-1}\Delta_{n,n_0}O.$$

Then, we have $(M_{n,n_0})^2 = V_n - V_{n_0-1}$. Let us denote by $|\cdot|_2$ the usual Euclidean norm on $\mathbb{R}^d$ and $\|A\| := \sup_{|z|_2=1} |Az|_2$ for any $(d,d)$-matrix $A$. Let us recall that if $A$ is a nonnegative symmetric matrix, then $\|A\|$ is equal to the maximal eigenvalue of $A$. We have

$$|\mathbb{E}[\phi(\sqrt{n}Y)] - \mathbb{E}[\phi(T_{n_0-1,n})]| = |\mathbb{E}[\phi(\sqrt{n}Y) - \phi(M_{n,n_0}Y)]|$$
$$\leq L_\phi \mathbb{E}[|(\sqrt{n}I_d - M_{n,n_0})Y|_\infty]$$
$$\leq L_\phi \mathbb{E}[|(\sqrt{n}I_d - M_{n,n_0})Y|_2]$$
$$\leq L_\phi \|\sqrt{n}I_d - M_{n,n_0}\| \mathbb{E}[|Y|_2]$$
$$\leq L_\phi \sqrt{\|nI_d - (V_n - V_{n_0-1})\|} \mathbb{E}[|Y|_2]$$
$$\leq L_\phi \sqrt{\|nI_d - V_n\| + \|V_{n_0-1}\|} \mathbb{E}[|Y|_2].$$

Indeed we have

$$\|\sqrt{n}I_d - M_{n,n_0}\| = \max_{\lambda \in \mathrm{Sp}(M_{n,n_0})} |\sqrt{n} - \lambda| = \max_{\mu \in \mathrm{Sp}(V_n - V_{n_0-1})} |\sqrt{n} - \sqrt{\mu}|$$

$$\leq \max_{\mu \in \mathrm{Sp}(V_n - V_{n_0-1})} \sqrt{|n - \mu|} = \sqrt{\|nI_d - (V_n - V_{n_0-1})\|},$$

where we denote by $\mathrm{Sp}(A)$ the set of eigenvalues of the square matrix $A$. Therefore, since the sequence of matrices $(V_m - m \cdot I_d)_m$ is bounded, we have

$$|\mathbb{E}[\phi(S_n)] - \mathbb{E}[\phi(\sqrt{n}Y)]| \leq K'(1 + A\psi(\varepsilon) + 2\varepsilon)L_\phi,$$

with $K' := K_0 + \mathbb{E}[|Y|_2] \sup_{m \geq n_0} \sqrt{\|mI_d - V_m\| + \|V_{n_0-1}\|}$. Let us denote by $\varepsilon_A$ the unique real number $\varepsilon_A \in [1, +\infty[$ such that $1 + A\psi(\varepsilon_A) = \varepsilon_A$.

According to the preceding, for any integer $n \geq 9n_0$ and any real number $A \geq M$, we have

$$(\mathcal{P}_n(A)) \quad \Longrightarrow \quad (\mathcal{P}_{n+1}(3K'\varepsilon_A)).$$

Let us show that there exists a real number $A_0 \geq M$ such that, for any $A \geq A_0$, we have $3K'\varepsilon_A \leq A$. The function $A \mapsto \varepsilon_A$ is increasing. If we had $M_1 := \sup_A \varepsilon_A < +\infty$, we would have $m_1 := \inf_A \psi(\varepsilon_A) > 0$ and therefore, for any $A$, $M_1 \geq \varepsilon_A \geq A\psi(\varepsilon_A) \geq Am_1$, which is impossible. Therefore, we have $\lim_{A \to +\infty} \varepsilon_A = +\infty$. Hence, there exists $A_0 \geq M$ such that for any $A \geq A_0$, we have $3K'(1 + A\psi(\varepsilon_A)) \leq A$ and therefore $3K'\varepsilon_A \leq A$.

Hence, for any real number $A \geq \max(A_0, A_1)$, we have $(\mathcal{P}_{9n_0}(A))$ and, for any integer $n \geq 9n_0$, $(\mathcal{P}_n(A)) \Rightarrow (\mathcal{P}_{n+1}(A))$, from which we deduce Theorem 1.1.

Now, we have to prove Proposition A.1.

Let an integer $n \geq 9n_0$ be given. Let $(Y_k)_{k \geq n_0}$ be a sequence of independent random variables defined on $(\Omega, \mathcal{F}, \nu)$ independent of $(X_k)_{k \geq 0}$ such that $Y_k$ is a Gaussian random variable with expectation 0 and covariance matrix $v_k$. Let $Y$ be a Gaussian random variable with expectation 0 and covariance matrix $I_d$, defined on $(\Omega, \mathcal{F}, \nu)$, independent of $((Y_k)_{k \geq n_0}, (X_k)_{k \geq 0})$ [this is always possible in some extension of $(\Omega, \mathcal{F}, \nu)$].

NOTATION A.2.    Let $k$ be an integer such that $n_0 \leq k \leq n$. We define $\Delta_k(f) = \mathbb{E}[f(S_{k-1} + X_k)] - \mathbb{E}[f(S_{k-1} + Y_k)]$.

For any function $\phi \in \mathbf{Lip}(\mathbb{R}^d, \mathbb{R})$, for any real number $\varepsilon \geq 1$ and any $x \in \mathbb{R}^d$, we define

$$f_{\phi,k,n,\varepsilon}(x) := \mathbb{E}\left[\phi\left(x + \sum_{i=k+1}^n Y_i + \varepsilon Y\right)\right],$$

with convention $\sum_{i=k+1}^n Y_i = 0$ if $k = n$.

Let us notice that we have

$$\mathbb{E}[\phi(S_n + \varepsilon Y)] - \mathbb{E}\left[\phi\left(S_{n_0-1} + \sum_{i=n_0}^{n} Y_i + \varepsilon Y\right)\right] = \sum_{k=n_0}^{n} \Delta_k(f_{\phi,k,n,\varepsilon}).$$

We will use Taylor expansions for functions $h \colon \mathbb{R}^d \to \mathbb{R}$. We will use the following notation.

NOTATION A.3.   Let $k$ be an integer such that $k \geq 1$.
If $h \colon \mathbb{R}^d \to \mathbb{R}$ is $k$-times differentiable, for any $x \in \mathbb{R}^d$, we denote by $D^k h(x)$ the point of $\mathbb{R}^{\{1,\dots,d\}^k}$ given by

$$D^k h(x) := \left(\frac{\partial^k}{\partial x_{i_1} \cdots \partial x_{i_k}} h(x)\right)_{i_1,\dots,i_k=1,\dots,d}.$$

We denote by $|\cdot|_\infty$ the supremum norm on $\mathbb{R}^{\{1,\dots,d\}^k}$.
For any $A^{(1)}, \dots, A^{(k)}$ in $\mathbb{R}^d$, we denote by $A^{(1)} \otimes \cdots \otimes A^{(k)}$ the point of $\mathbb{R}^{\{1,\dots,d\}^k}$ given by

$$A^{(1)} \otimes \cdots \otimes A^{(k)} = \left(\prod_{j=1}^{k} A_{i_j}^{(j)}\right)_{i_1,\dots,i_k=1,\dots,d}$$

[if $A^{(l)} = (A_1^{(l)}, \dots, A_d^{(l)})$].
For any integer $j$ satisfying $1 \leq j \leq k$, for any $A \in \mathbb{R}^{\{1,\dots,d\}^k}$ and any $B \in \mathbb{R}^{\{1,\dots,d\}^j}$, we denote by $A * B$ the point of $\mathbb{R}^{\{1,\dots,d\}^{k-j}}$ given by

$$A * B := \left(\sum_{n_1,\dots,n_j=1}^{d} A_{n_1,\dots,n_j,i_1,\dots,i_{k-j}} B_{n_1,\dots,n_j}\right)_{i_1,\dots,i_{k-j}=1,\dots,d}.$$

For any $A \colon \Omega \to \mathbb{R}^{\{1,\dots,d\}^k}$ and any $B \colon \Omega \to \mathbb{R}^{\{1,\dots,d\}^k}$, we define (when it is well defined):

$$\mathbb{E}[A] := (\mathbb{E}[A_{i_1,\dots,i_k}])_{i_1,\dots,i_k=1,\dots,d},$$

$$\|A\|_\infty := \||A|_\infty\|_\infty,$$

$$\mathrm{Cov}(A, B) := \sum_{i_1,\dots,i_k=1}^{d} \mathrm{Cov}(A_{i_1,\dots,i_k}, B_{i_1,\dots,i_k}).$$

Let us notice that if $j = k$, $*$ corresponds to the usual scalar product on $\mathbb{R}^{\{1,\dots,d\}^k}$.
On the other hand, let us notice that the $k$-linear form on $\mathbb{R}^d$ associated to $D^k h(x)$ is

$$(A_1, \dots, A_k) \mapsto D^k h(x) * (A_1 \otimes \cdots \otimes A_k)$$

and that, for any $j = 1, \ldots, k$, we have

$$(D^k h(x)) * (A_1 \otimes \cdots \otimes A_k) = ((D^k h(x)) * (A_1 \otimes \cdots \otimes A_j)) * (A_{j+1} \otimes \cdots \otimes A_k).$$

We have $\Delta_k(f) = \Delta_{1,k}(f) - \Delta_{2,k}(f)$, with

$$\Delta_{1,k}(f) := \mathbb{E}[f(S_{k-1} + X_k)] - \mathbb{E}[f(S_{k-1})] - \tfrac{1}{2}\mathbb{E}[D^2 f(S_{k-1})] * v_k$$

and

$$\Delta_{2,k}(f) := \mathbb{E}[f(S_{k-1} + Y_k)] - \mathbb{E}[f(S_{k-1})] - \tfrac{1}{2}\mathbb{E}[D^2 f(S_{k-1})] * v_k.$$

**A.1. Estimations for small $k$.**

LEMMA A.1.1 (Adaptation of Lemma 6 of [37]). *Let* $f : \mathbb{R}^d \to \mathbb{R}$ *be a function in* $C^4$. *We have*

$$|\Delta_k(f)| \le d^4 (\|D^3 f\|_\infty + \|D^4 f\|_\infty) \Bigg( M^3 + 15 C^2 M^2 (r+1) \sum_{p=0}^{k-1} (1+p) \varphi_{p,0}$$
$$+ 3CM \sum_{p=r+1}^{k-2} \sum_{l=1,\ldots,k-1 \,:\, (r+1)l \le p} \varphi_{p,l} \Bigg).$$

PROOF. Since $Y_k$ is a Gaussian random variable independent of $S_{k-1}$, with expectation 0 and covariance matrix $v_k$, we have

$$|\Delta_{2,k}(f)| = \left| \mathbb{E}\left[ f(S_{k-1} + Y_k) - f(S_{k-1}) - \tfrac{1}{2} D^2 f(S_{k-1}) * (Y_k \otimes Y_k) \right] \right|$$

$$= \frac{1}{2} \left| \int_0^1 (1-t)^2 \mathbb{E}[D^3 f(S_{k-1} + tY_k) * Y_k^{\otimes 3}]\, dt \right|$$

$$= \frac{1}{2} \left| \int_0^1 (1-t)^2 \mathbb{E}[g(tY_k) * Y_k^{\otimes 3}]\, dt \right|$$

$$\qquad\qquad\qquad \text{with } g(u) = \mathbb{E}[D^3 f(S_{k-1} + u)]$$

$$\le \frac{d^3}{2} \int_0^1 (1-t)^2 \sup_{a \in \mathbb{R}^d} |\mathbb{E}[D^3 f(S_{k-1} + a)]|_\infty \mathbb{E}[|Y_k^{\otimes 3}|_\infty]\, dt$$

$$\le \frac{d^3}{6} \sup_{a \in \mathbb{R}^d} |\mathbb{E}[D^3 f(S_{k-1} + a)]|_\infty \mathbb{E}[|Y_k^{\otimes 3}|_\infty].$$

We have $\mathbb{E}[|Y_k^{\otimes 3}|_\infty] \le \frac{4d|v_k|_\infty^{3/2}}{\sqrt{2\pi}}$. Moreover, according to hypothesis (1), we have $|v_k|_\infty \le 2C(M+1) \sum_{p=0}^{k-1} \varphi_{p,0}$. According to Hölder inequality and to the fact that $\varphi_{p,0} \leq 1$, we can show (cf. [37], page 264) that we have $|v_k|_\infty^{3/2} \leq (4C(M+1))^{3/2} \frac{\pi}{\sqrt{6}} \sum_{p=0}^{k-1} (1+p)\varphi_{p,0}$. Hence, we have

$$
\begin{aligned}
|\Delta_{2,k}(f)| \leq d^4 \frac{8\sqrt{\pi}}{3\sqrt{3}} &\sup_{a \in \mathbb{R}} |\mathbb{E}[D^3 f(S_{k-1} + a)]|_\infty \\
&\times (C(M+1))^{3/2} \sum_{p=0}^{k-1} (1+p)\varphi_{p,0}.
\end{aligned}
\tag{7}
$$

Let us now control $\Delta_{1,k}(f)$. Since we have $v_k = \mathbb{E}[X_k^{\otimes 2}] + \sum_{i=1}^{k-1} (\mathbb{E}[X_i \otimes X_k] + \mathbb{E}[X_k \otimes X_i])$, we have

$$
\begin{aligned}
\Delta_{1,k}(f) = \mathbb{E}[D^1 f(S_{k-1}) * X_k] &+ \tfrac{1}{2} \mathrm{Cov}(D^2 f(S_{k-1}), X_k^{\otimes 2}) \\
&- \mathbb{E}[D^2 f(S_{k-1})] * \sum_{i=1}^{k-1} \mathbb{E}[X_i \otimes X_k] \\
&+ \mathbb{E}[\tfrac{1}{6} D^3 f(S_{k-1} + \theta_k X_k) * X_k^{\otimes 3}],
\end{aligned}
$$

where $\theta_k$ is a random variable with values in $[0;1]$. We have

$$
\|\tfrac{1}{6} D^3 f(S_{k-1} + \theta_k X_k) * X_k^{\otimes 3}\|_\infty \leq \tfrac{1}{6} d^3 \|D^3 f\|_\infty M^3.
\tag{8}
$$

On the other hand, according to (1), we have

$$
\begin{aligned}
|\mathrm{Cov}(D^2 f(S_{k-1}), X_k^{\otimes 2})| &\leq \sum_{i=1}^{k-1} |\mathrm{Cov}(D^2 f(S_i) - D^2 f(S_{i-1}), X_k^{\otimes 2})| \\
&\leq 3 d^3 C M \|D^3 f\|_\infty \sum_{p=1}^{k-1} \varphi_{p,0}.
\end{aligned}
\tag{9}
$$

We have

$$
\begin{aligned}
D^1 f(S_{k-1}) &= D^1 f(0) + \sum_{i=1}^{k-1} (D^1 f(S_i) - D^1 f(S_{i-1})) \\
&= D^1 f(0) + \sum_{i=1}^{k-1} \Big( D^2 f(S_{i-1}) * X_i \\
&\qquad\qquad + \int_0^1 (1-t) D^3 f(S_{i-1} + t X_i) * X_i^{\otimes 2} \, dt \Big)
\end{aligned}
$$

and so

$$
\mathbb{E}[D^1 f(S_{k-1}) * X_k] - \mathbb{E}[D^2 f(S_{k-1})] * \sum_{i=1}^{k-1} \mathbb{E}[X_i \otimes X_k]
$$

$$(10) \quad = \sum_{i=1}^{k-1} \mathrm{Cov}(D^2 f(S_{i-1}), X_i \otimes X_k)$$

$$+ \sum_{i=1}^{k-1} \mathbb{E}[D^2 f(S_{i-1}) - D^2 f(S_{k-1})] * \mathbb{E}[X_i \otimes X_k]$$

$$+ \sum_{i=1}^{k-1} \mathrm{Cov}\left(\int_0^1 (1-t) D^3 f(S_{i-1} + tX_i) * (X_i^{\otimes 2}) \, dt, X_k\right),$$

since $\mathbb{E}[D^1 f(0) * X_k] = D^1 f(0) * \mathbb{E}[X_k] = 0$. According to (1), we have

$$\sum_{i=1}^{k-1} |\mathbb{E}[D^2 f(S_{i-1}) - D^2 f(S_{k-1})] * \mathbb{E}[X_i \otimes X_k]|$$

$$(11) \quad \leq d^3 \|D^3 f\|_\infty M 2 C M \sum_{i=1}^{k-1} (k-i) \varphi_{k-i,0}$$

$$\leq 2 d^3 \|D^3 f\|_\infty C M^2 \sum_{p=1}^{k-1} p \varphi_{p,0}$$

and

$$\sum_{i=1}^{k-1} \left| \mathrm{Cov}\left(\int_0^1 (1-t) D^3 f(S_{i-1} + tX_i) * (X_i^{\otimes 2}) \, dt, X_k\right) \right|$$

$$(12) \quad \leq \sum_{i=1}^{k-1} d^3 C(2\|D^3 f\|_\infty M^2 + \|D^4 f\|_\infty M^2) \varphi_{k-i,0}$$

$$\leq 2 C d^3 (\|D^3 f\|_\infty + \|D^4 f\|_\infty) M^2 \sum_{p=1}^{k-1} \varphi_{p,0}.$$

For any integer $i = 1, \ldots, k-1$, we write $j = j_i := \max(0, (r+2)i - (r+1)k)$. According to (1), we have

$$|\mathrm{Cov}((D^2 f(S_{i-1}) - D^2 f(S_j)) * X_i, X_k)|$$

$$\leq \sum_{m=j+1}^{i-1} |\mathrm{Cov}((D^2 f(S_m) - D^2 f(S_{m-1})) * X_i, X_k)|$$

$$\leq 3 d^3 C \|D^3 f\|_\infty M^2 (i-j-1) \varphi_{k-i,0}$$

and

$$|\mathbb{E}[D^2 f(S_{i-1}) - D^2 f(S_j)] * \mathbb{E}[X_i \otimes X_k]|$$

$$\leq d^3 \|D^3 f\|_\infty (i-j-1) M 2 C M \varphi_{k-i,0}.$$

Hence, we have

$$
(13) \quad \begin{aligned}
&\sum_{i=1}^{k-1} |\mathrm{Cov}(D^2 f(S_{i-1}) - D^2 f(S_j), X_i \otimes X_k)| \\
&\leq 5d^3 C M^2 \|D^3 f\|_\infty (r+1) \sum_{p=1}^{k-1} p\varphi_{p,0}.
\end{aligned}
$$

If $(r+2)i - (r+1)k \leq 0$, then $j = 0$ and so $\mathrm{Cov}(D^2 f(S_j), X_i \otimes X_k) = 0$. Hence we have

$$
(14) \quad \begin{aligned}
&\sum_{i=1}^{k-1} |\mathrm{Cov}(D^2 f(S_j), X_i \otimes X_k)| \\
&\leq \sum_{i=1,\ldots,k-1\,:\,(r+1)k < (r+2)i} \sum_{l=1}^{j} |\mathrm{Cov}(D^2 f(S_l) - D^2 f(S_{l-1}), X_i \otimes X_k)| \\
&\leq \sum_{i=1,\ldots,k-1\,:\,(r+1)k < (r+2)i} d^3 3 C M \|D^3 f\|_\infty \sum_{l=1}^{(r+2)i - (r+1)k} \varphi_{i-l,k-i} \\
&\leq 3d^3 C M \|D^3 f\|_\infty \sum_{p=r+1}^{k-2} \sum_{j=1,\ldots,k-1\,:\,(r+1)j \leq p} \varphi_{p,j}. \qquad \square
\end{aligned}
$$

LEMMA A.1.2 (Adaptation of Lemma 5 of [37]). *For any $\phi \in \mathbf{Lip}(\mathbb{R}^d, \mathbb{R})$, for any integer $k = n_0, \ldots, n$ and any real number $\varepsilon \geq 1$, the function $f_{\phi,k,n,\varepsilon}$ is $C^\infty$ and, for any integer $i \geq 1$, we have*

$$
\|D^i f_{\phi,k,n,\varepsilon}\|_\infty \leq \frac{C_i}{(n - k + \varepsilon^2)^{(i-1)/2}} L_\phi,
$$

*with $C_i := d^{i+1} 3 \cdot 2^{i-1} \int_{\mathbb{R}^d} |z|_\infty |D^i h(z)|_\infty \, dz$, where $h$ is the density function of the Gaussian law with expectation 0 and covariance matrix $I_d$.*

PROOF. Let us denote by $\Gamma_{n,k,\varepsilon^2}$ the positive symmetric matrix such that $\Gamma_{n,k,\varepsilon^2}^2 = V_n - V_k + \varepsilon^2 I_d$. For any $x \in \mathbb{R}^d$, we have

$$
\begin{aligned}
f_{\phi,k,n,\varepsilon}(x) &= \mathbb{E}[\phi(x + Y_{k+1} + \cdots + Y_n + \varepsilon Y)] \\
&= \frac{1}{\sqrt{\det(V_n - V_k + \varepsilon^2 I_d)}} \int_{\mathbb{R}^d} \phi(u) h(\Gamma_{n,k,\varepsilon^2}^{-1}(u - x)) \, du.
\end{aligned}
$$

Let an integer $i \geq 1$ be given. Let $k_1, \ldots, k_i$ in $\{1, \ldots, d\}$ be given. Let us denote, for any matrix $A$, the $j$th column vector of $A$ by $[A]_j$. For any $x \in \mathbb{R}^d$,

we have

$$(-1)^i (D^i f_{\phi,k,n,\varepsilon}(x))_{k_1,\dots,k_i}$$

$$= \frac{1}{\sqrt{\det(V_n - V_k + \varepsilon^2 I_d)}}$$

$$\times \int_{\mathbb{R}^d} \phi(u) D^i h(\Gamma_{n,k,\varepsilon^2}^{-1}(u-x)) * ([\Gamma_{n,k,\varepsilon^2}^{-1}]_{k_1} \otimes \cdots \otimes [\Gamma_{n,k,\varepsilon^2}^{-1}]_{k_i}) \, du$$

$$= \int_{\mathbb{R}^d} \phi(x + \Gamma_{n,k,\varepsilon^2} \cdot z) D^i h(z) * ([\Gamma_{n,k,\varepsilon^2}^{-1}]_{k_1} \otimes \cdots \otimes [\Gamma_{n,k,\varepsilon^2}^{-1}]_{k_i}) \, dz$$

$$= \int_{\mathbb{R}^d} (\phi(x + \Gamma_{n,k,\varepsilon^2} \cdot z) - \phi(x))$$

$$\times D^i h(z) * ([\Gamma_{n,k,\varepsilon^2}^{-1}]_{k_1} \otimes \cdots \otimes [\Gamma_{n,k,\varepsilon^2}^{-1}]_{k_i}) \, dz,$$

since $\int_{\mathbb{R}^d} D^i h(z) \, dz = 0$. Let us denote by $\lambda_-$ and $\lambda_+$, respectively, the smallest and biggest eigenvalues of $V_n - V_k + \varepsilon^2 I_d$. Since $\Gamma_{n,k,\varepsilon^2}$ is diagonalizable in an orthonormal basis, we have

$$|\Gamma_{n,k,\varepsilon^2} \cdot z|_\infty \le d\lambda_+^{1/2}|z|_\infty.$$

Moreover, for any $j = 1, \dots, d$, we have $|[\Gamma_{n,k,\varepsilon^2}^{-1}]_j|_\infty \le \lambda_-^{-1/2}$. Therefore we have

$$|D^i f_{\phi,k,n,\varepsilon^2}(x)|_\infty \le L_\phi d \frac{\lambda_+^{1/2} d^i}{\lambda_-^{i/2}} \int_{\mathbb{R}^d} |z|_\infty |D^i h(z)| \, dz.$$

Since $k \ge n_0$, we have

$$\tfrac{1}{2}(n-k) + \varepsilon^2 \le \lambda_- \le \lambda_+ \le \tfrac{3}{2}(n-k) + \varepsilon^2$$

(according to the fact that two invertible symmetric matrices are diagonal in a same basis). $\square$

PROPOSITION A.1.3.    *For any* $\phi \in \mathbf{Lip}(\mathbb{R}^d, \mathbb{R})$, *for any real number* $\varepsilon \ge 1$, *we have*

$$\sum_{k=n_0}^{n-\lfloor n/3 \rfloor - 1} |\Delta_k(f_{\phi,k,n,\varepsilon})| \le d^4 (C_3 + C_4) \ln(3)$$

$$\times \left( M^3 + 15 C^2 M^2 (r+1) \sum_{p \ge 0} (1+p)\varphi_{p,0} \right.$$

$$\left. + 3CM \sum_{p \ge r+1} \sum_{l=1}^{\lfloor p/(r+1) \rfloor} \varphi_{p,l} \right) L_\phi.$$

PROOF. According to Lemmas A.1.1 and A.1.2, we have

$$\sum_{k=n_0}^{n-\lfloor n/3 \rfloor-1} |\Delta_k(f_{\phi,k,n,\varepsilon})|$$

$$\leq d^4(C_3+C_4)$$

$$\times \sum_{m=\lfloor n/3 \rfloor+1}^{n-1} \frac{1}{m+1}\Bigg(M^3+15C^2M^2(r+1)\sum_{p\geq 0}(1+p)\varphi_{p,0}$$

$$+3CM\sum_{p\geq r+1}\sum_{l=1}^{\lfloor p/(r+1) \rfloor}\varphi_{p,l}\Bigg)L_\phi. \quad \square$$

## A.2. Estimations for big $k$.

LEMMA A.2.1 (Analogous to Lemma 7 of [37]). *For any real number* $A \geq M$, *for any integer* $n \geq 9n_0$, *if property* $(\mathcal{P}_n(A))$ *is satisfied, then, for any* $\phi \in \mathbf{Lip}(\mathbb{R}^d, \mathbb{R})$, *for any real number* $\varepsilon \geq 1$, *for any integer* $k \in [n - \lfloor\frac{n}{3}\rfloor; n]$, *for any integer* $l \in [\frac{n}{3}; k]$, *for any integer* $i \geq 1$, *we have*

$$(15) \quad \sup_{a\in\mathbb{R}^d} |\mathbb{E}[D^i f_{\phi,k,n,\varepsilon}(S_l+a)]|_\infty \leq K_i\left(\frac{A}{(n-k+\varepsilon^2)^{i/2}}+\frac{1}{n^{(i-1)/2}}\right)L_\phi,$$

*with* $K_i := 2^{i/2}d^i\int_{\mathbb{R}^d}|D^i h(z)|_\infty\,dz + C_i\sqrt{3}^{i-1}$.

PROOF. Let $N_l$ be a Gaussian random variable with null expectation and covariance matrix $l \cdot I_d$. Let $\Gamma_{n,k,\varepsilon^2}$ be as in the proof of Lemma A.1.2. First, let us notice that, if $(\mathcal{P}_n(A))$ is satisfied, we have

$$|\mathbb{E}[D^i f_{\phi,k,n,\varepsilon}(S_l+a)] - \mathbb{E}[D^i f_{\phi,k,n,\varepsilon}(N_l+a)]|_\infty$$

$$\leq \lambda_-^{-i/2}d^i\int_{\mathbb{R}^d}|\mathbb{E}[\phi(S_l+a+\Gamma_{n,k,\varepsilon^2}\cdot z)]$$

$$-\mathbb{E}[\phi(N_l+a+\Gamma_{n,k,\varepsilon^2}\cdot z)]|\cdot|D^i h(z)|_\infty\,dz$$

$$\leq \frac{2^{i/2}d^i}{(n-k+\varepsilon^2)^{i/2}}\int_{\mathbb{R}^d}AL_\phi|D^i h(z)|_\infty\,dz,$$

since $y \mapsto \phi(y+a+\Gamma_{n,k,\varepsilon^2}z)$ is Lipschitz continuous with Lipschitz constant bounded by $L_\phi$. On the other hand, we have

$$\mathbb{E}[D^i f_{\phi,k,n,\varepsilon}(N_l+a)] = D^i(a \mapsto (\mathbb{E}[f_{\phi,k,n,\varepsilon}(N_l+a)]))$$

$$= D^i(a \mapsto (\mathbb{E}[\phi(a+\Gamma_{n,k,\varepsilon^2+l}\cdot N)])),$$

where $N$ is a Gaussian random variable with null expectation and covariance matrix $I_d$. As in the proof of Lemma A.1.2, we get

$$|\mathbb{E}[D^i f_{\phi,k,n,\varepsilon}(N_l + a)]|_\infty \leq \frac{C_i}{\sqrt{n-k+l+\varepsilon^2}^{i-1}} L_\phi \leq \frac{C_i \sqrt{3}^{i-1}}{\sqrt{n}^{i-1}} L_\phi. \qquad \square$$

According to estimation (7) and since $k - 1 \geq \frac{n}{3}$, we have:

LEMMA A.2.2. *Let an integer $n \geq 9n_0$ and a real number $A \geq M$ be given. If property $(\mathcal{P}_n(A))$ is satisfied, then, for any $\phi \in \mathbf{Lip}(\mathbb{R}^d, \mathbb{R})$, for any integer $k = n - \lfloor \frac{n}{3} \rfloor, \ldots, n$, for any real number $\varepsilon \geq 1$, we have*

$$\begin{aligned}
(16) \quad |\Delta_{2,k}(f_{\phi,k,n,\varepsilon})| \leq &\, L_\phi \cdot d^4 \frac{16\sqrt{2\pi}}{3\sqrt{3}} K_3 \left( \frac{A}{(n-k+\varepsilon^2)^{3/2}} + \frac{1}{n} \right) (CM)^{3/2} \\
&\times \sum_{p=0}^{k-1} (1+p)\varphi_{p,0}.
\end{aligned}$$

LEMMA A.2.3. *There exists a real number $\tilde{K}$ (only depending on $d$, $C$, $M$ and $r$) such that, for any $\phi \in \mathbf{Lip}(\mathbb{R}^d, \mathbb{R})$, for any integer $n \geq 9n_0$ and any real number $A \geq M$, if property $(\mathcal{P}_n(A))$ is satisfied, then, for any integer $k = n - \lfloor \frac{n}{3} \rfloor, \ldots, n$, we have*

$$\begin{aligned}
(17) \quad |\Delta_{1,k}(f_{\phi,k,n,\varepsilon})| \leq &\, L_\phi \tilde{K} \left( \frac{A\alpha_{\sqrt{n-k}}}{(n-k+\varepsilon^2)^{3/2}} + \frac{\beta_{\sqrt{n-k}}}{n} \right. \\
&\left. + \frac{A\gamma_{\sqrt{n-k}}}{n-k+\varepsilon^2} + \frac{A\delta_{\sqrt{n-k}}}{\sqrt{n-k+\varepsilon^2}} + \varphi_{\lfloor \sqrt{n-k} \rfloor + 1, 0} \right),
\end{aligned}$$

*with $\alpha_m := 1 + \sum_{p=1}^{\lfloor m \rfloor} p\zeta_p$, $\beta_m := 1 + \sum_{p=1}^{\lfloor m \rfloor} \zeta_p$, $\gamma_m := \sum_{p=\lceil m/(r+2)^2 \rceil}^{\lfloor m \rfloor} \zeta_p$ and $\delta_m := \sum_{p=\lfloor m/(r+2) \rfloor + 1}^{+\infty} \frac{\zeta_p}{p}$, with $\zeta_p := p \max_{j=0,\ldots,\lfloor p/(r+1) \rfloor} \varphi_{p,j}$.*

PROOF. To simplify notation, we will write $f_k$ instead of $f_{\phi,k,n,\varepsilon}$. We have

$$\begin{aligned}
\Delta_{1,k}(f_k) = &\, \mathbb{E}[D^1 f_k(S_{k-1}) * X_k] + \tfrac{1}{2}\mathbb{E}[D^2 f_k(S_{k-1}) * (X_k^{\otimes 2} - v_k)] \\
&+ \tfrac{1}{6}\mathbb{E}[D^3 f_k(S_{k-1}) * X_k^{\otimes 3}] \\
&+ \tfrac{1}{6}\mathbb{E}\left[ \int_0^1 (1-t)^3 D^4 f_k(S_{k-1} + tX_k) * (X_k^{\otimes 4}) \, dt \right].
\end{aligned}$$

Since we have $v_k = \mathbb{E}[X_k^{\otimes 2}] + \sum_{i=1,\dots,k-1}(\mathbb{E}[X_k \otimes X_i] + \mathbb{E}[X_i \otimes X_k])$ and since the matrix $\mathbb{E}[D^2 f_k(S_{k-1})]$ is symmetric, we have

$$
\begin{aligned}
\Delta_{1,k}(f_k) = {} & \mathbb{E}[D^1 f_k(S_{k-1}) * X_k] - \sum_{i=1,\dots,k-1} \mathbb{E}[D^2 f_k(S_{k-1})] * \mathbb{E}[X_k \otimes X_i] \\
& + \tfrac{1}{2}\operatorname{Cov}(D^2 f_k(S_{k-1}), X_k^{\otimes 2}) \\
& + \tfrac{1}{6}\mathbb{E}[D^3 f_k(S_{k-1}) * X_k^{\otimes 3}] \\
& + \tfrac{1}{6}\mathbb{E}\left[\int_0^1 (1-t)^3 D^4 f_k(S_{k-1} + tX_k) * X_k^{\otimes 4}\, dt\right].
\end{aligned}
$$

For any integer $\ell \in \{0,\dots,k-1\}$, we have

$$
\begin{aligned}
D^1 f_k & (S_{k-1}) - D^1 f_k(S_{k-\ell-1}) \\
& = \sum_{j=k-\ell}^{k-1} \left(D^1 f_k(S_j) - D^1 f_k(S_{j-1})\right) \\
& = \sum_{j=k-\ell}^{k-1} D^2 f_k(S_{j-1}) * X_j \\
& \quad + \tfrac{1}{2}\sum_{j=k-\ell}^{k-1} D^3 f_k(S_{j-1}) * X_j^{\otimes 2} \\
& \quad + \tfrac{1}{2}\sum_{j=k-\ell}^{k-1} \int_0^1 (1-t)^2 D^4 f_k(S_{j-1} + tX_j) * X_j^{\otimes 3}\, dt.
\end{aligned}
$$

Therefore we have

$$
\begin{aligned}
\Delta_{1,k}(f_k) = {} & \mathbb{E}\left[\tfrac{1}{6}\int_0^1 (1-t)^3 D^4 f_k(S_{k-1} + tX_k) * X_k^{\otimes 4}\, dt\right] \\
& + \tfrac{1}{6}\mathbb{E}[D^3 f_k(S_{k-1}) * X_k^{\otimes 3}] \\
& + \mathbb{E}[D^1 f_k(S_{k-\ell-1}) * X_k] \\
& - \sum_{i=1}^{k-\ell-1} \mathbb{E}[D^2 f_k(S_{k-1})] * \mathbb{E}[X_k \otimes X_i] \\
& + \tfrac{1}{2}\operatorname{Cov}(D^2 f_k(S_{k-1}), X_k^{\otimes 2}) \\
& + \tfrac{1}{2}\sum_{j=k-\ell}^{k-1} \mathbb{E}[D^3 f_k(S_{j-1}) * (X_j^{\otimes 2} \otimes X_k)]
\end{aligned}
$$

$$+ \frac{1}{2} \sum_{j=k-\ell}^{k-1} \mathbb{E}\left[\int_0^1 (1-t)^2 D^4 f_k(S_{j-1} + tX_j) * (X_k \otimes X_j^{\otimes 3}) \, dt\right]$$

$$- \sum_{j=k-\ell}^{k-1} \mathbb{E}[D^2 f_k(S_{k-1}) - D^2 f_k(S_{j-1})] * \mathbb{E}[X_k \otimes X_j]$$

$$+ \sum_{j=k-\ell}^{k-1} \mathrm{Cov}(D^2 f_k(S_{j-1}), X_j \otimes X_k).$$

In the following, we take $\ell = \lfloor \sqrt{n-k} \rfloor$, the integer part of $\sqrt{n-k}$. Since $n \geq 9n_0 \geq 9$ and $\frac{2n}{3} \leq k \leq n$, we have $0 \leq \ell \leq k-1$.

A.2.1. *Control of* $\mathbb{E}[\frac{1}{6} \int_0^1 (1-t)^3 D^4 f_k(S_{k-1} + tX_k) X_k^{\otimes 4} \, dt]$. According to Lemma A.1.2, we have

$$\left| \mathbb{E}\left[ \frac{1}{6} \int_0^1 (1-t)^3 D^4 f_k(S_{k-1} + tX_k) * X_k^{\otimes 4} \, dt \right] \right|$$

(18)
$$\leq d^4 \frac{1}{24} \|D^4 f_k\|_\infty M^4$$

$$\leq d^4 \frac{M^4 C_4}{24} \frac{\alpha_{\sqrt{n-k}}}{(n-k+\varepsilon^2)^{3/2}} L_\phi.$$

A.2.2. *Control of* $\frac{1}{6}\mathbb{E}[D^3 f_k(S_{k-1}) * X_k^{\otimes 3}]$. We have $D^3 f_k(S_{k-1}) = D^3 f_k(S_{k-\ell-1}) + \sum_{j=1}^{\ell} (D^3 f_k(S_{k-j}) - D^3 f_k(S_{k-j-1}))$. According to (1) and to Lemma A.1.2, we have

(19)     $|\mathrm{Cov}(D^3 f_k(S_{k-\ell-1}), X_k^{\otimes 3})| \leq d^3 L_\phi C(C_3 + C_4)\varphi_{\lfloor \sqrt{n-k} \rfloor + 1, 0}$

and

$$\sum_{j=1}^{\ell} |\mathrm{Cov}((D^3 f_k(S_{k-j}) - D^3 f_k(S_{k-j-1})), X_k^{\otimes 3})|$$

(20)
$$\leq d^3 \sum_{j=1}^{\ell} \frac{3C_4 M}{(n-k+\varepsilon^2)^{3/2}} \varphi_{j,0} L_\phi$$

$$\leq d^3 L_\phi 3 C_4 M \frac{A \alpha_{\sqrt{n-k}}}{(n-k+\varepsilon^2)^{3/2}}.$$

Moreover, according to Lemma A.2.1, we have

$$|\mathbb{E}[D^3 f_k(S_{k-1})] * \mathbb{E}[X_k^{\otimes 3}]|$$

(21)
$$\leq d^3 K_3 \left( \frac{A}{(n-k+\varepsilon^2)^{3/2}} + \frac{1}{n} \right) M^3 L_\phi$$

$$\leq d^3 K_3 M^3 \left( \frac{A\alpha_{\sqrt{n-k}}}{(n-k+\varepsilon^2)^{3/2}} + \frac{\beta_{\sqrt{n-k}}}{n} \right) L_\phi.$$

A.2.3. *Control of* $\mathbb{E}[D^1 f_k(S_{k-\ell-1}) * X_k]$. According to (1) and to Lemma A.1.2, we have

$$|\mathbb{E}[D^1 f_k(S_{k-\ell-1}) * X_k]| = |\mathrm{Cov}(D^1 f_k(S_{k-\ell-1}), X_k)|$$

$$\leq d(C_1 + C_2) \varphi_{\lfloor\sqrt{n-k}\rfloor+1,0} L_\phi. \tag{22}$$

A.2.4. *Control of* $\sum_{j=\ell+1}^{k-1} |\mathbb{E}[D^2 f_k(S_{k-1})] * \mathbb{E}[X_{k-j} \otimes X_k]|$. We have

$$\sum_{j=\ell+1}^{k-1} |\mathbb{E}[D^2 f_k(S_{k-1})] * \mathbb{E}[X_{k-j} \otimes X_k]|$$

$$\leq d^2 \sum_{j=\ell+1}^{k-1} |\mathbb{E}[D^2 f_k(S_{k-1})]|_\infty |\mathbb{E}[X_{k-j} \otimes X_k]|_\infty$$

$$\leq d^2 \sum_{j=\ell+1}^{k-1} \frac{C_2}{\sqrt{n-k+\varepsilon^2}} L_\phi 2CM \varphi_{j,0} \tag{23}$$

$$\leq d^2 2CMC_2 \frac{\delta_{\sqrt{n-k}}}{\sqrt{n-k+\varepsilon^2}} L_\phi.$$

A.2.5. *Control of* $|\mathrm{Cov}(D^2 f_k(S_{k-1}), X_k^{\otimes 2})|$. We have

$$D^2 f_k(S_{k-1}) = D^2 f_k(S_{k-\ell-1}) + \sum_{j=1}^{\ell} (D^2 f_k(S_{k-j}) - D^2 f_k(S_{k-j-1}))$$

$$= D^2 f_k(S_{k-\ell-1}) + \sum_{j=1}^{\ell} D^3 f_k(S_{k-j-1}) * X_{k-j}$$

$$+ \sum_{j=1}^{\ell} \int_0^1 (1-t) D^4 f_k(S_{k-j-1} + tX_{k-j}) * X_{k-j}^{\otimes 2} dt.$$

Hence, we have

$$\mathrm{Cov}(D^2 f_k(S_{k-1}), X_k^{\otimes 2})$$

$$= \mathrm{Cov}(D^2 f_k(S_{k-\ell-1}), X_k^{\otimes 2}) + \sum_{j=1}^{\ell} \mathrm{Cov}(D^3 f_k(S_{k-\ell-1}) * X_{k-j}, X_k^{\otimes 2})$$

$$+ \sum_{j=1}^{\ell} \sum_{m=j+1}^{\ell} \mathrm{Cov}((D^3 f_k(S_{k-m}) - D^3 f_k(S_{k-m-1})) * X_{k-j}, X_k^{\otimes 2})$$

$$+ \sum_{j=1}^{\ell} \int_0^1 (1-t) \, \mathrm{Cov}(D^4 f_k(S_{k-j-1} + tX_{k-j}) * X_{k-j}^{\otimes 2}, X_k^{\otimes 2}) \, dt.$$

1. First, according to (1) and to Lemma A.1.2, we have

$$
\begin{aligned}
|\mathrm{Cov}(D^2 f_k(S_{k-\ell-1}), X_k^{\otimes 2})| &\le d^2 C \frac{C_2 + C_3}{\sqrt{n-k+\varepsilon^2}} L_\phi \varphi_{l+1,0} \\
&\le d^2 C (C_2 + C_3) \varphi_{\lfloor \sqrt{n-k} \rfloor + 1, 0} L_\phi.
\end{aligned}
\tag{24}
$$

2. Control of $\sum_{j=1}^{\ell} \mathrm{Cov}(D^3 f_k(S_{k-\ell-1}) * X_{k-j}, X_k^{\otimes 2})$.

(a) For any integer $j = 1, \dots, \ell$ satisfying $\frac{\sqrt{n-k}}{r+2} < j \le \ell$, we have $\sqrt{n-k} < (r+2)j$ and, according to (1) and to Lemma A.1.2, we have

$$|\mathrm{Cov}(D^3 f_k(S_{k-\ell-1}) * X_{k-j}, X_k^{\otimes 2})| \le d^3 \frac{2C(C_3 + C_4)M}{n-k+\varepsilon^2} L_\phi \varphi_{j,0}.$$

Hence

$$
\begin{aligned}
\sum_{j=\lfloor \sqrt{n-k}/(r+2) \rfloor + 1}^{\ell} &|\mathrm{Cov}(D^3 f_k(S_{k-\ell-1}) * X_{k-j}, X_k^{\otimes 2})| \\
&\le d^3 \frac{2C(C_3 + C_4)M}{n-k+\varepsilon^2} \gamma_{\sqrt{n-k}} L_\phi.
\end{aligned}
\tag{25}
$$

(b) For any integer $j = 1, \dots, \ell$ satisfying $j \le \frac{\sqrt{n-k}}{r+2}$, according to (1) and to Lemma A.1.2, we have

$$|\mathrm{Cov}(D^3 f_k(S_{k-\ell-1}), X_{k-j} \otimes X_k^{\otimes 2})| \le d^3 \frac{C(C_3 + C_4)}{n-k+\varepsilon^2} L_\phi \varphi_{\ell+1-j,j}$$

and

$$|\mathrm{Cov}(D^3 f_k(S_{k-\ell-1}), X_{k-j}) * \mathbb{E}[X_k^{\otimes 2}]| \le d^3 \frac{CM^2(C_3 + C_4)}{n-k+\varepsilon^2} L_\phi \varphi_{\ell+1-j,0}.$$

Moreover, according to (1) and to Lemma A.2.1, since $k - \ell - 1 \ge \frac{n}{3}$, we have

$$
\begin{aligned}
|\mathbb{E}[D^3 f_k(S_{k-\ell-1})] &* \mathbb{E}[X_{k-j} \otimes X_k^{\otimes 2}]| \\
&\le d^3 K_3 \left( \frac{A}{(n-k+\varepsilon^2)^{3/2}} + \frac{1}{n} \right) L_\phi C 2M \varphi_{j,0},
\end{aligned}
$$

from which we get

$$\sum_{j=1}^{\lfloor \sqrt{n-k}/(r+2) \rfloor} |\mathrm{Cov}(D^3 f_k(S_{k-\ell-1}) * X_{k-j}, X_k^{\otimes 2})|$$

$$(26) \qquad \leq d^3 2CM^2(C_3+C_4)\frac{\gamma_{\sqrt{n-k}}}{n-k+\varepsilon^2}L_\phi$$

$$+ 2CMK_3\frac{A\alpha_{\sqrt{n-k}}}{(n-k+\varepsilon^2)^{3/2}}L_\phi + K_3 2CM\frac{\beta_{\sqrt{n-k}}}{n}L_\phi.$$

3. Control of $\sum_{j=1}^{\ell}\sum_{m=j+1}^{\ell}\mathrm{Cov}((D^3f_k(S_{k-m})-D^3f_k(S_{k-m-1}))*X_{k-j},X_k^{\otimes 2})$.

(a) For any integers $j$ and $m$ satisfying $1 \leq j < m \leq \ell$ and $m \leq (r+2)j$, according to (1) and to Lemma A.1.2, we have

$$|\mathrm{Cov}((D^3f_k(S_{k-m})-D^3f_k(S_{k-m-1}))*X_{k-j},X_k^{\otimes 2})|$$

$$\leq d^4\frac{3CC_4M^2}{(n-k+\varepsilon^2)^{3/2}}L_\phi\varphi_{j,0}.$$

Hence, we have

$$\sum_{j=1}^{\ell}\sum_{m=j+1}^{\min(\ell,(r+2)j)}|\mathrm{Cov}((D^3f_k(S_{k-m})-D^3f_k(S_{k-m-1}))*X_{k-j},X_k^{\otimes 2})|$$

$$(27)$$

$$\leq d^4\frac{3CC_4M^2(r+1)\alpha_{\sqrt{n-k}}}{(n-k+\varepsilon^2)^{3/2}}L_\phi.$$

(b) We will use the following formula: $\mathrm{Cov}(A*B,C) = \mathrm{Cov}(A,B \otimes C) - \mathrm{Cov}(A,B)*\mathbb{E}[C] + \mathbb{E}[A]*\mathbb{E}[B \otimes (C-\mathbb{E}[C])]$. For any integers $j$ and $m$ satisfying $1 \leq j \leq (r+2)j+1 \leq m \leq \ell$, according to (1) and to Lemma A.1.2, we have

$$|\mathrm{Cov}(D^3f_k(S_{k-m})-D^3f_k(S_{k-m-1}),X_{k-j}\otimes X_k^{\otimes 2})|$$

$$\leq d^4C\frac{3C_4M}{(n-k+\varepsilon^2)^{3/2}}L_\phi\varphi_{m-j,j}$$

and

$$|\mathrm{Cov}(D^3f_k(S_{k-m})-D^3f_k(S_{k-m-1}),X_{k-j})*\mathbb{E}[X_k^{\otimes 2}]|$$

$$\leq d^4C\frac{3C_4M^3}{(n-k+\varepsilon^2)^{3/2}}L_\phi\varphi_{m-j,0},$$

from which we get

$$\sum_{j=1}^{\ell}\sum_{(r+2)j+1\leq m\leq \ell}|\mathrm{Cov}(D^3f_k(S_{k-m})-D^3f_k(S_{k-m-1}),X_{k-j}\otimes X_k^{\otimes 2})|$$

$$+\sum_{j=1}^{\ell}\sum_{(r+2)j+1\leq m\leq \ell}|\mathrm{Cov}(D^3f_k(S_{k-m})$$

(28)
$$- D^3 f_k(S_{k-m-1}), X_{k-j}) * \mathbb{E}[X_k^{\otimes 2}]|$$

$$\leq d^4 \frac{6CC_4 M^3}{(r+1)(n-k+\varepsilon^2)^{3/2}} \alpha_{\sqrt{n-k}} L_\phi.$$

Indeed, we have

$$\sum_{j=1}^{\ell} \sum_{(r+2)j+1 \leq m \leq \ell} \varphi_{m-j,j} \leq \sum_{p=r+2}^{\ell} \sum_{j=1}^{\lfloor p/(r+1) \rfloor} \varphi_{p,j} \leq \frac{\alpha_{\sqrt{n-k}}}{r+1}$$

and

$$\sum_{j=1}^{\ell} \sum_{(r+2)j+1 \leq m \leq \ell} \varphi_{m-j,0} \leq \frac{\alpha_{\sqrt{n-k}}}{r+1}.$$

Finally, according to (1) and to Lemma A.2.1 (since $k - \ell - 1 \geq \frac{n}{3}$), we have

$$\sum_{j=1}^{\ell} \left| \sum_{m=\min(\ell,(r+2)j)+1}^{\ell} \mathbb{E}[D^3 f_k(S_{k-m}) - D^3 f_k(S_{k-m-1})] \right.$$

$$* \left. \mathbb{E}[X_{k-j} \otimes (X_k^{\otimes 2} - \mathbb{E}[X_k^{\otimes 2}])] \right|$$

$$= \sum_{j=1}^{\ell} |\mathbb{E}[D^3 f_k(S_{k-\min(\ell,(r+2)j)-1}) - D^3 f_k(S_{k-\ell-1})]$$

(29)
$$* \mathbb{E}[X_{k-j} \otimes (X_k^{\otimes 2} - \mathbb{E}[X_k^{\otimes 2}])]|$$

$$= \sum_{j=1}^{\lfloor \ell/(r+2) \rfloor} |\mathbb{E}[D^3 f_k(S_{k-(r+2)j-1}) - D^3 f_k(S_{k-\ell-1})]$$

$$* \mathbb{E}[X_{k-j} \otimes (X_k^{\otimes 2} - \mathbb{E}[X_k^{\otimes 2}])]|$$

$$\leq d^3 \sum_{j=1}^{\lfloor \ell/(r+2) \rfloor} 2K_3 \left( \frac{A}{(n-k+\varepsilon^2)^{3/2}} + \frac{1}{n} \right) L_\phi 2CM \varphi_{j,0}$$

$$\leq d^3 4CK_3 M \left( \frac{A\alpha_{\sqrt{n-k}}}{(n-k+\varepsilon^2)^{3/2}} + \frac{\beta_{\sqrt{n-k}}}{n} \right) L_\phi.$$

4. We have

$$\sum_{j=1}^{\ell} \int_0^1 (1-t) |\text{Cov}(D^4 f_k(S_{k-j-1} + tX_{k-j}) * X_{k-j}^{\otimes 2}, X_k^{\otimes 2})| \, dt$$

$$(30) \qquad \leq d^4 \sum_{j=1}^{\ell} \int_0^1 (1-t) \frac{3C(C_4+C_5)M^2}{(n-k+\varepsilon^2)^{3/2}} L_\phi \varphi_{j,0} \, dt$$

$$\leq d^4 \frac{3C(C_4+C_5)M^2}{(n-k+\varepsilon^2)^{3/2}} \alpha_{\sqrt{n-k}} L_\phi.$$

A.2.6. *Control of* $\sum_{j=1}^{\ell} \mathbb{E}[D^3 f_k(S_{k-j-1}) * (X_{k-j}^{\otimes 2} \otimes X_k)]$. Let us notice that we have

$$\sum_{j=1}^{\ell} \mathbb{E}[D^3 f_k(S_{k-j-1}) * X_{k-j}^{\otimes 2} \otimes X_k] = \sum_{j=1}^{\ell} \mathrm{Cov}(D^3 f_k(S_{k-j-1}) * X_{k-j}^{\otimes 2}, X_k).$$

We control this quantity as we did for

$$\sum_{j=1}^{\ell} \mathrm{Cov}(D^3 f_k(S_{k-j-1}) * X_{k-j}, X_k^{\otimes 2})$$

in the previous section (we get analogous estimations).

A.2.7. *Control of* $\sum_{j=1}^{\ell} \int_0^1 (1-t)^2 \mathbb{E}[D^4 f_k(S_{k-j-1} + tX_{k-j}) * (X_{k-j}^{\otimes 3} \otimes X_k)] \, dt$. We control quantity $\sum_{j=1}^{\ell} \int_0^1 (1-t)^2 \, \mathrm{Cov}(D^4 f_k(S_{k-j-1} + tX_{k-j}) * X_{k-j}^{\otimes 3}, X_k) \, dt$ as we did for $\sum_{j=1}^{\ell} \int_0^1 (1-t) \, \mathrm{Cov}(D^4 f_k(S_{k-j-1} + tX_{k-j}) * X_{k-j}^{\otimes 2}, X_k^{\otimes 2}) \, dt$. We obtain estimations analogous to (30).

A.2.8. *Control of* $\sum_{j=1}^{\ell} \mathbb{E}[D^2 f_k(S_{k-1}) - D^2 f_k(S_{k-j-1})] * \mathbb{E}[X_{k-j} \otimes X_k]$. For any integer $j = 1, \ldots, \ell$, we have

$$D^2 f_k(S_{k-1}) - D^2 f_k(S_{k-j-1})$$

$$= \sum_{m=1}^{j} (D^2 f_k(S_{k-m}) - D^2 f_k(S_{k-m-1}))$$

$$= \sum_{m=1}^{j} \Big( D^3 f_k(S_{k-m-1}) * X_{k-m}$$

$$+ \int_0^1 (1-t) D^4 f(S_{k-m-1} + tX_{k-m}) * X_{k-m}^{\otimes 2} \, dt \Big).$$

Moreover we have

$$D^3 f_k(S_{k-m-1}) = D^3 f_k(S_{k-\ell-1}) + \sum_{p=m+1}^{\ell} (D^3 f_k(S_{k-p}) - D^3 f_k(S_{k-p-1})).$$

We get

$$\sum_{j=1}^{\ell} \mathbb{E}[D^2 f_k(S_{k-1}) - D^2 f_k(S_{k-j-1})] * \mathbb{E}[X_{k-j} \otimes X_k]$$

$$= \sum_{j=1}^{\ell} \sum_{m=1}^{j} \mathbb{E}\left[\int_0^1 (1-t) D^4 f_k(S_{k-m-1} + t X_{k-m}) * X_{k-m}^{\otimes 2} \, dt\right]$$

$$* \mathbb{E}[X_{k-j} \otimes X_k]$$

$$+ \sum_{j=1}^{\ell} \sum_{m=1}^{j} \mathbb{E}[D^3 f_k(S_{k-\ell-1}) * X_{k-m}] * \mathbb{E}[X_{k-j} \otimes X_k]$$

$$+ \sum_{j=1}^{\ell} \sum_{m=1}^{j} \sum_{p=m+1}^{\ell} \mathbb{E}[D^3(f_k(S_{k-p}) - D^3 f_k(S_{k-p-1})) * X_{k-m}]$$

$$* \mathbb{E}[X_{k-j} \otimes X_k].$$

1. According to (1) and to Lemma A.1.2, we have

$$\sum_{j=1}^{\ell} \sum_{m=1}^{j} \left| \mathbb{E}\left[\int_0^1 (1-t) D^4 f_k(S_{k-m-1} + t X_{k-m}) * X_{k-m}^{\otimes 2} \, dt\right]\right.$$

$$\left. * \mathbb{E}[X_{k-j} \otimes X_k] \right|$$

(31)

$$\leq d^4 \sum_{j=1}^{\ell} \sum_{m=1}^{j} \frac{C_4}{(n-k+\varepsilon^2)^{3/2}} L_\phi M^2 C 2M \varphi_{j,0}$$

$$\leq d^4 \frac{2CM^3 C_4}{(n-k+\varepsilon^2)^{3/2}} \alpha_{\sqrt{n-k}} L_\phi.$$

2. We have $\sum_{j=1}^{\ell} \sum_{m=1}^{j} \mathbb{E}[D^3 f_k(S_{k-\ell-1}) * X_{k-m}] * \mathbb{E}[X_{k-j} \otimes X_k] = A_1 + A_2$ with

$$A_1 = \sum_{j=1}^{\ell} \sum_{m=1}^{j} \mathbb{E}[D^3 f_k(S_{k-\ell-1}) * X_{k-m}] * \mathbb{E}[X_{k-j} \otimes X_k] \mathbf{1}_{\{\ell \geq (r+2)m\}},$$

$$A_2 = \sum_{j=1}^{\ell} \sum_{m=1}^{j} \mathbb{E}[D^3 f_k(S_{k-\ell-1}) * X_{k-m}] * \mathbb{E}[X_{k-j} \otimes X_k] \mathbf{1}_{\{\ell < (r+2)m\}}.$$

According to (1) and to Lemma A.1.2, we have

$$|A_1| \leq \sum_{j=1}^{\ell} \sum_{m=1}^{j} d^3 C M^2 \frac{C_3 + C_4}{n - k + \varepsilon^2} L_\phi \varphi_{\ell+1-m,0} \mathbf{1}_{\{\ell \geq (r+2)m\}}$$

$$\leq d^3 C M^2 \frac{C_3 + C_4}{n - k + \varepsilon^2} L_\phi \sum_{m=1}^{\lfloor \sqrt{n-k}/(r+2) \rfloor} (\ell + 1 - m) \varphi_{\ell+1-m,0}$$

$$(32) \qquad \leq d^3 C M^2 \frac{C_3 + C_4}{n - k + \varepsilon^2} L_\phi \sum_{p=\lceil (r+1)\sqrt{n-k}/(r+2) \rceil}^{\ell} p \varphi_{p,0}$$

$$(33) \qquad \leq d^3 C M^2 \frac{C_3 + C_4}{n - k + \varepsilon^2} L_\phi \gamma_{\sqrt{n-k}}.$$

[We use the fact that if $\ell \geq (r+2)m$, then we have $m \leq \frac{\sqrt{n-k}}{r+2}$ and $\ell + 1 - m \geq \frac{(r+1)\sqrt{n-k}}{r+2} \geq (r+1)m$ and $m \leq \frac{\ell+1-m}{r+1}$.]

On the other hand, we have

$$|A_2| \leq d^3 \sum_{j=1}^{\ell} \sum_{m=1}^{j} \frac{2 C M^2 C_3}{n - k + \varepsilon^2} \varphi_{j,0} \mathbf{1}_{\{\ell < (r+2)m\}} L_\phi$$

$$(34) \qquad \leq d^3 \frac{2 C M^2 C_3}{n - k + \varepsilon^2} \sum_{j=\lceil \sqrt{n-k}/(r+2) \rceil}^{\ell} j \varphi_{j,0} L_\phi$$

$$\leq d^3 \frac{2 C M^2 C_3}{n - k + \varepsilon^2} \gamma_{\sqrt{n-k}} L_\phi.$$

3. We have

$$\sum_{j=1}^{\ell} \sum_{m=1}^{j} \sum_{p=m+1}^{(r+2)j} |\mathbb{E}[(D^3 f_k(S_{k-p}) - D^3 f_k(S_{k-p-1})) * X_{k-m}]$$

$$* \mathbb{E}[X_{k-j} \otimes X_k]|$$

$$(35)$$

$$\leq \sum_{j=1}^{\ell} \sum_{m=1}^{j} \sum_{p=m+1}^{(r+2)j} \frac{d^4 C_4 M^2}{(n - k + \varepsilon^2)^{3/2}} L_\phi 2 C M \varphi_{j,0}$$

$$\leq \frac{2 d^4 C C_4 M^3}{(n - k + \varepsilon^2)^{3/2}} (r+2) \alpha_{\sqrt{n-k}} L_\phi$$

and

$$\sum_{j=1}^{\ell} \sum_{m=1}^{j} \sum_{p=(r+2)j+1}^{\ell} |\mathbb{E}[(D^3 f_k(S_{k-p}) - D^3 f_k(S_{k-p-1})) * X_{k-m}]$$

$$* \mathbb{E}[X_{k-j} \otimes X_k]|$$

$$(36) \qquad \leq d^4 \sum_{j=1}^{\ell} \sum_{m=1}^{j} \sum_{p=(r+2)j+1}^{\ell} \frac{3C_4 M^3}{(n-k+\varepsilon^2)^{3/2}} L_\phi \varphi_{p-m,0}$$

$$\leq d^4 \frac{3C_4 M^3}{(n-k+\varepsilon^2)^{3/2}} \frac{1}{(r+1)^2} \sum_{p=1}^{\ell} p^2 \varphi_{p,0} L_\phi$$

$$\leq d^4 \frac{3C_4 M^3}{(n-k+\varepsilon^2)^{3/2}} \frac{1}{(r+1)^2} \alpha_{\sqrt{n-k}} L_\phi.$$

A.2.9. *Control of* $\sum_{j=1}^{\ell} \mathrm{Cov}(D^2 f_k(S_{k-j-1}), X_{k-j} \otimes X_k)$. For any integer $j = 1, \ldots, \ell$, we have

$$D^2 f_k(S_{k-j-1}) = D^2 f_k(S_{k-\ell-1}) + \sum_{m=j+1}^{\ell} (D^2 f_k(S_{k-m}) - D^2 f_k(S_{k-m-1}))$$

and, for any integer $m = j+1, \ldots, \ell$:

$$D^2 f_k(S_{k-m}) - D^2 f_k(S_{k-m-1})$$

$$= D^3 f_k(S_{k-m-1}) * X_{k-m} + \int_0^1 (1-t) D^4 f(S_{k-m-1} + tX_{k-m}) * X_{k-m}^{\otimes 2} \, dt$$

$$= D^3 f_k(S_{k-\ell-1}) * X_{k-m} + \sum_{p=m+1}^{\ell} (D^3 f_k(S_{k-p}) - D^3 f_k(S_{k-p-1})) * X_{k-m}$$

$$+ \int_0^1 (1-t) D^4 f(S_{k-m-1} + tX_{k-m}) * X_{k-m}^{\otimes 2} \, dt.$$

Therefore, we have

$$\sum_{j=1}^{\ell} \mathrm{Cov}(D^2 f_k(S_{k-j-1}), X_{k-j} \otimes X_k)$$

$$= \sum_{j=1}^{\ell} \mathrm{Cov}(D^2 f_k(S_{k-\ell-1}), X_{k-j} \otimes X_k)$$

$$+ \sum_{j=1}^{\ell} \sum_{m=j+1}^{\ell} \mathrm{Cov}\Big( \int_0^1 (1-t) D^4 f_k(S_{k-m-1} + tX_{k-m}) X_{k-m}^{\otimes 2} \, dt,$$

$$X_{k-j} \otimes X_k \Big)$$

$$+ \sum_{j=1}^{\ell} \sum_{m=j+1}^{\ell} \mathrm{Cov}(D^3 f_k(S_{k-\ell-1}) * X_{k-m}, X_{k-j} \otimes X_k) \mathbf{1}_{\{\ell < (r+2)m\}}$$

$$+ \sum_{j=1}^{\ell} \sum_{m=j+1}^{\ell} \mathrm{Cov}(D^3 f_k(S_{k-\ell-1}) * X_{k-m}, X_{k-j} \otimes X_k) \mathbf{1}_{\{\ell \geq (r+2)m\}}$$

$$+ \sum_{j=1}^{\ell} \sum_{m=j+1}^{\ell} \sum_{p=m+1}^{\min(\ell,(r+2)m)} \mathrm{Cov}((D^3 f_k(S_{k-p}) - D^3 f_k(S_{k-p-1}))$$
$$* X_{k-m}, X_{k-j} \otimes X_k)$$

$$+ \sum_{j=1}^{\ell} \sum_{m=j+1}^{\ell} \mathrm{Cov}((D^3 f_k(S_{k-(r+2)m-1}) - D^3 f_k(S_{k-\ell-1}))$$
$$* X_{k-m}, X_{k-j} \otimes X_k) \mathbf{1}_{\{\ell \geq (r+2)m+1\}}.$$

Let us control each term of the right-hand member of this identity.

1. [Control of $\sum_{j=1}^{\ell} \mathrm{Cov}(D^2 f_k(S_{k-\ell-1}), X_{k-j} \otimes X_k)$.]
According to (1) and to Lemma A.1.2, we have

$$\sum_{j=\lfloor \ell/(r+2) \rfloor + 1}^{\ell} |\mathrm{Cov}(D^2 f_k(S_{k-\ell-1}), X_{k-j} \otimes X_k)|$$

$$\leq \sum_{j=\lfloor \ell/(r+2) \rfloor + 1}^{\ell} (|\mathrm{Cov}(D^2 f_k(S_{k-\ell-1}) * X_{k-j}, X_k)|$$
$$+ |\mathbb{E}[D^2 f_k(S_{k-\ell-1})] * \mathbb{E}[X_{k-j} \otimes X_k]|)$$

$$\leq \sum_{j=\lfloor \ell/(r+2) \rfloor + 1}^{\ell} d^2 \frac{4C(C_2 + C_3)M}{\sqrt{n-k+\varepsilon^2}} L_\phi \varphi_{j,0}$$

and

$$\sum_{j=1}^{\lfloor \ell/(r+2) \rfloor} |\mathrm{Cov}(D^2 f_k(S_{k-\ell-1}), X_{k-j} \otimes X_k)|$$

$$\leq \sum_{j=1}^{\lfloor \ell/(r+2) \rfloor} d^2 \frac{C(C_2 + C_3)}{\sqrt{n-k+\varepsilon^2}} L_\phi \varphi_{\ell+1-j,j}.$$

Let us notice that if $j \leq \lfloor \frac{\ell}{r+2} \rfloor$, then we have $(r+2)j \leq \ell$ and so $j = \frac{(r+1)j}{r+1} \leq \frac{\ell-j+1}{r+1}$ and $\ell + 1 - j \geq \ell - \lfloor \frac{\ell}{r+2} \rfloor + 1 \geq \frac{r+1}{r+2}\sqrt{n-k} \geq \frac{\sqrt{n-k}}{r+2}$. Therefore, we have

$$\sum_{j=1}^{\ell} |\mathrm{Cov}(D^2 f_k(S_{k-\ell-1}), X_{k-j} \otimes X_k)|$$

(37)
$$\leq d^2 \frac{5C(C_2+C_3)M}{\sqrt{n-k+\varepsilon^2}}\delta_{\sqrt{n-k}}L_\phi.$$

2. [Control of $\sum_{j=1}^{\ell}\sum_{m=j+1}^{\ell}\mathrm{Cov}(\int_0^1(1-t)D^4f_k(S_{k-m-1}+tX_{k-m})*X_{k-m}^{\otimes 2}\,dt, X_{k-j}\otimes X_k).$]

Let $j$ and $m$ be two integers satisfying $1\leq j\leq j+1\leq m\leq\ell$.

(a) If $m\leq(r+2)j$, then we have

$$\left|\mathrm{Cov}\left(\int_0^1(1-t)D^4f_k(S_{k-m-1}+tX_{k-m})*(X_{k-m}^{\otimes 2}\otimes X_{k-j})\,dt, X_k\right)\right|$$

$$+\left|\mathbb{E}\left[\int_0^1(1-t)D^4f_k(S_{k-m-1}+tX_{k-m})*X_{k-m}^{\otimes 2}\,dt\right]\right.$$

$$\left.*\mathbb{E}[X_{k-j}\otimes X_k]\right|$$

$$\leq d^4C\frac{3(C_4+C_5)M^3}{(n-k+\varepsilon^2)^{3/2}}L_\phi\varphi_{j,0}.$$

(b) If $m>(r+2)j$, then we have

$$\left|\mathrm{Cov}\left(\int_0^1(1-t)D^4f_k(S_{k-m-1}+tX_{k-m})*X_{k-m}^{\otimes 2}\,dt, X_{k-j}\otimes X_k\right)\right|$$

$$\leq d^4C\frac{2(C_4+C_5)M^2}{(n-k+\varepsilon^2)^{3/2}}L_\phi\varphi_{m-j,j}.$$

Therefore we have

$$\sum_{j=1}^{\ell}\sum_{m=j+1}^{\ell}\left|\mathrm{Cov}\left(\int_0^1(1-t)D^4f_k(S_{k-m-1}+tX_{k-m})*X_{k-m}^{\otimes 2}\,dt, X_{k-j}\otimes X_k\right)\right|$$

(38)
$$\leq\frac{\tilde{K}_0}{(n-k+\varepsilon^2)^{3/2}}\left(\sum_{p=1}^{\lfloor\sqrt{n-k}\rfloor}p\varphi_{p,0}+\sum_{p=1}^{\lfloor\sqrt{n-k}\rfloor}\sum_{j=1}^{\lfloor p/(r+1)\rfloor}\varphi_{p,j}\right)L_\phi$$

$$\leq\frac{2\tilde{K}_0\alpha_{\sqrt{n-k}}}{(n-k+\varepsilon^2)^{3/2}}L_\phi,$$

for some $\tilde{K}_0$ only depending on $d$, $C$, $C_4$, $C_5$, $M$ and $r$.

3. We have

$$\sum_{j=1}^{\ell}\sum_{m=j+1}^{\ell}|\mathrm{Cov}(D^3f_k(S_{k-\ell-1})*X_{k-m}, X_{k-j}\otimes X_k)\mathbf{1}_{\{\ell<(r+2)m\}}|$$

$$\leq d^3\frac{\tilde{K}_1}{n-k+\varepsilon^2}\left(\sum_{p=\lceil\sqrt{n-k}/(r+2)^2\rceil}^{\lfloor\sqrt{n-k}\rfloor}p\varphi_{p,0}\right.$$

(39)
$$+ \sum_{p=\lceil (r+1)\sqrt{n-k}/(r+2)^2 \rceil}^{\lfloor \sqrt{n-k} \rfloor} \sum_{j=1}^{\lfloor p/(r+1) \rfloor} \varphi_{p,j} \Bigg) L_\phi$$

$$\leq d^3 \frac{2\tilde{K}_1 \gamma_{\sqrt{n-k}}}{n-k+\varepsilon^2} L_\phi,$$

for some $\tilde{K}_1$ only depending on $d$, $C$, $M$ and $r$. Indeed, let $j$ and $m$ be two integers satisfying $1 \leq j \leq j+1 \leq m \leq \ell$.

(a) If $\ell < (r+2)m$ and $m \leq (r+2)j$, then we have

$$|\mathrm{Cov}(D^3 f_k(S_{k-\ell-1}) * (X_{k-m} \otimes X_{k-j}), X_k)| \leq d^3 C \frac{2(C_3+C_4)}{n-k+\varepsilon^2} L_\phi M^2 \varphi_{j,0}$$

and

$$|\mathbb{E}[D^3 f_k(S_{k-\ell-1}) * X_{k-m}] * \mathbb{E}[X_{k-j} \otimes X_k]| \leq d^3 \frac{C_3 M}{n-k+\varepsilon^2} C 2 M \varphi_{j,0} L_\phi.$$

(b) If $\ell < (r+2)m$ and $m > (r+2)j$, then we have

$$|\mathrm{Cov}(D^3 f_k(S_{k-\ell-1}) * X_{k-m}, X_{k-j} \otimes X_k)| \leq d^3 C \frac{2(C_3+C_4)M}{n-k+\varepsilon^2} L_\phi \varphi_{m-j,j}.$$

4. We have

$$\sum_{j=1}^{\ell} \sum_{m=j+1}^{\ell} \mathrm{Cov}(D^3 f_k(S_{k-\ell-1}) * X_{k-m}, X_{k-j} \otimes X_k) \mathbf{1}_{\{\ell \geq (r+2)m\}}$$

$$= \sum_{j=1}^{\ell} \sum_{m=j+1}^{\ell} \mathbb{E}[D^3 f_k(S_{k-\ell-1})] * \mathbb{E}[X_{k-m} \otimes X_{k-j} \otimes X_k] \mathbf{1}_{\{\ell \geq (r+2)m\}}$$

$$+ \sum_{j=1}^{\ell} \sum_{m=j+1}^{\ell} \mathrm{Cov}(D^3 f_k(S_{k-\ell-1}), X_{k-m} \otimes X_{k-j} \otimes X_k) \mathbf{1}_{\{\ell \geq (r+2)m\}}$$

$$- \sum_{j=1}^{\ell} \sum_{m=j+1}^{\ell} \mathbb{E}[D^3 f_k(S_{k-\ell-1}) * X_{k-m}] * \mathbb{E}[X_{k-j} \otimes X_k] \mathbf{1}_{\{\ell \geq (r+2)m\}}.$$

(a) We have

$$\sum_{j=1}^{\ell} \sum_{m=j+1}^{\ell} |\mathbb{E}[D^3 f_k(S_{k-\ell-1})] * \mathbb{E}[X_{k-m} \otimes X_{k-j} \otimes X_k]|$$

$$\leq d^3 K_3 \left( \frac{A}{(n-k+\varepsilon^2)^{3/2}} + \frac{1}{n} \right) L_\phi$$

$$\times \sum_{j=1}^{\ell} \sum_{m=j+1}^{\ell} |\mathbb{E}[X_{k-m} \otimes X_{k-j} \otimes X_k]|_\infty.$$

Since we have

$$\sum_{j=1}^{\ell}\sum_{m=j+1}^{(r+2)j}|\mathbb{E}[X_{k-m}\otimes X_{k-j}\otimes X_k]|_\infty \leq \sum_{j=1}^{\ell}\sum_{m=j+1}^{(r+2)j}2CM^2\varphi_{j,0}$$

$$\leq 2CM^2\sum_{j=1}^{\ell}(r+1)j\varphi_{j,0}$$

and

$$\sum_{j=1}^{\ell}\sum_{m=(r+2)j+1}^{\ell}|\mathbb{E}[X_{k-m}\otimes X_{k-j}\otimes X_k]|_\infty \leq \sum_{j=1}^{\ell}\sum_{m=(r+2)j+1}^{\ell}2CM\varphi_{m-j,j}$$

$$\leq 2CM\sum_{p=1}^{\ell}\sum_{j=1}^{\lfloor p/(r+1)\rfloor}\varphi_{p,j},$$

we get

(40)
$$\sum_{j=1}^{\ell}\sum_{m=j+1}^{\ell}|\mathbb{E}[D^3 f_k(S_{k-\ell-1})]*\mathbb{E}[X_{k-m}\otimes X_{k-j}\otimes X_k]|$$

$$\leq d^3 4(r+1)CM^2 K_3\left(\frac{A\alpha\sqrt{n-k}}{(n-k+\varepsilon^2)^{3/2}}+\frac{\beta\sqrt{n-k}}{n}\right)L_\phi.$$

(b) For the two other terms, we write

$$|\mathrm{Cov}(D^3 f_k(S_{k-\ell-1}),X_{k-m}\otimes X_{k-j}\otimes X_k)|\leq d^3 C\frac{C_3+C_4}{n-k+\varepsilon^2}L_\phi\varphi_{\ell+1-m,m}$$

and

$$|\mathbb{E}[D^3 f_k(S_{k-\ell-1})*X_{k-m}]*\mathbb{E}[X_{k-j}\otimes X_k]|$$

$$\leq d^3 CM^2\frac{C_3+C_4}{n-k+\varepsilon^2}L_\phi\varphi_{\ell+1-m,0},$$

according to (1) and to Lemma A.1.2. Let us notice that, if $\ell\geq(r+2)m$, then we have $m\leq\frac{\sqrt{n-k}}{r+2}$ and so $\ell+1-m\geq\frac{(r+1)\sqrt{n-k}}{r+2}\geq(r+1)m$ and $m\leq\frac{\ell+1-m}{r+1}$. Hence, we have

$$\sum_{j=1}^{\ell}\sum_{m=j+1}^{\ell}\varphi_{\ell+1-m,m}\mathbf{1}_{\{\ell\geq(r+2)m\}}$$

$$\leq\sum_{m=2}^{\ell}(m-1)\varphi_{\ell+1-m,m}\mathbf{1}_{\{\ell\geq(r+2)m\}}$$

$$\leq \sum_{m=2}^{\lfloor \sqrt{n-k}/(r+2) \rfloor} \frac{\ell+1-m}{r+1} \max_{j \leq (\ell+1-m)/(r+1)} \varphi_{\ell+1-m,j}$$

$$\leq \sum_{p=\lceil (r+1)\sqrt{n-k}/(r+2) \rceil}^{\ell} p \max_{j \leq p/(r+1)} \varphi_{p,j} \leq \gamma_{\sqrt{n-k}}.$$

Therefore, we have

$$(41) \quad \begin{aligned} &\sum_{j=1}^{\ell} \sum_{m=j+1}^{\ell} |\mathrm{Cov}(D^3 f_k(S_{k-\ell-1}), X_{k-m} \otimes X_{k-j} \otimes X_k)| \mathbf{1}_{\{\ell \geq (r+2)m\}} \\ &\leq d^3 C \frac{(C_3+C_4)\gamma_{\sqrt{n-k}}}{n-k+\varepsilon^2} L_\phi. \end{aligned}$$

In the same way, we get

$$(42) \quad \begin{aligned} &\sum_{j=1}^{\ell} \sum_{m=j+1}^{\ell} |\mathbb{E}[D^3 f_k(S_{k-\ell-1}) * X_{k-m}] * \mathbb{E}[X_{k-j} \otimes X_k]| \mathbf{1}_{\{\ell \geq (r+2)m\}} \\ &\leq d^3 C M^2 \frac{(C_3+C_4)\gamma_{\sqrt{n-k}}}{n-k+\varepsilon^2} L_\phi. \end{aligned}$$

5. We have

$$(43) \quad \begin{aligned} &\sum_{j=1}^{\ell} \sum_{m=j+1}^{\ell} \sum_{p=m+1}^{\min(\ell,(r+2)m)} |\mathrm{Cov}((D^3 f_k(S_{k-p}) - D^3 f_k(S_{k-p-1})) \\ &\hspace{5cm} * X_{k-m}, X_{k-j} \otimes X_k)| \\ &\leq \frac{\tilde{K}_3}{(n-k+\varepsilon^2)^{3/2}} L_\phi \left( \sum_{p=1}^{\lfloor \sqrt{n-k} \rfloor} p^2 \varphi_{p,0} + \sum_{p=1}^{\lfloor \sqrt{n-k} \rfloor} p \sum_{j=1}^{\lfloor p/(r+1) \rfloor} \varphi_{p,j} \right) \\ &\leq \frac{\tilde{K}'_3 \alpha_{\sqrt{n-k}}}{(n-k+\varepsilon^2)^{3/2}} L_\phi, \end{aligned}$$

for some $\tilde{K}_3$ and $\tilde{K}'_3$ only depending on $d$, $C$, $C_4$, $M$ and $r$. Indeed, let $j$, $m$ and $p$ be three integers satisfying $1 \leq j \leq j+1 \leq m \leq m+1 \leq p \leq \ell$.

(a) If we have $p \leq (r+2)m$ and $m \leq (r+2)j$, then we have

$$(44) \quad \begin{aligned} &|\mathrm{Cov}((D^3 f_k(S_{k-p}) - D^3 f_k(S_{k-p-1})) * (X_{k-m} \otimes X_{k-j}), X_k)| \\ &\leq d^4 C \frac{3C_4}{(n-k+\varepsilon^2)^{3/2}} L_\phi M^3 \varphi_{j,0} \end{aligned}$$

and

(45)
$$|\mathbb{E}[(D^3 f_k(S_{k-p}) - D^3 f_k(S_{k-p-1})) * X_{k-m}] * \mathbb{E}[X_{k-j} \otimes X_k]|$$
$$\leq d^4 C \frac{3C_4}{(n-k+\varepsilon^2)^{3/2}} L_\phi M^3 \varphi_{j,0}.$$

(b) On the other hand, if $p \leq (r+2)m$ and $m \geq (r+2)j + 1$, then we have

$$|\text{Cov}((D^3 f_k(S_{k-p}) - D^3 f_k(S_{k-p-1})) * X_{k-m}, X_{k-j} \otimes X_k)|$$
$$\leq d^4 C \frac{3C_4}{(n-k+\varepsilon^2)^{3/2}} L_\phi M^2 \varphi_{m-j,j}.$$

We conclude with the use of the following formulas:

$$\sum_{j=1}^{\ell} \sum_{m=j+1}^{(r+2)j} \sum_{p=m+1}^{(r+2)m} \varphi_{j,0} \leq (r+2)^3 \sum_{j=1}^{\ell} j^2 \varphi_{j,0}$$

and

$$\sum_{j=1}^{\ell} \sum_{m=(r+2)j+1}^{\ell} \sum_{p=m+1}^{(r+2)m} \varphi_{m-j,j} = \sum_{j=1}^{\ell} \sum_{m=(r+2)j+1}^{\ell} (r+1)m\varphi_{m-j,j}$$
$$\leq (r+2) \sum_{p=1}^{\ell} \sum_{j=1}^{\lfloor p/(r+1) \rfloor} p\varphi_{p,j}.$$

6. Now we control

$$\sum_{j=1}^{\ell} \sum_{m=j+1}^{\ell} \text{Cov}((D^3 f_k(S_{k-(r+2)m-1}) - D^3 f_k(S_{k-\ell-1}))$$
$$* X_{k-m}, X_{k-j} \otimes X_k)\mathbf{1}_{\{\ell \geq (r+2)m+1\}}.$$

If $\ell \geq (r+2)m + 1$, then we have

$$D^3 f_k(S_{k-(r+2)m-1}) - D^3 f_k(S_{k-\ell-1})$$
$$= \sum_{p=(r+2)m+1}^{\ell} (D^3 f_k(S_{k-p}) - D^3 f_k(S_{k-p-1})).$$

We will use the following formula:

$$\text{Cov}(A * B, C \otimes D) = \text{Cov}(A, B \otimes C \otimes D)$$
$$- \mathbb{E}[A * B] * \mathbb{E}[C \otimes D] + \mathbb{E}[A] * \mathbb{E}[B \otimes C \otimes D].$$

(a) If $\ell \geq p \geq (r+2)m+1$, then we have

$$|\mathrm{Cov}(D^3 f_k(S_{k-p}) - D^3 f_k(S_{k-p-1}), X_{k-m} \otimes X_{k-j} \otimes X_k)|$$

$$\leq d^4 C \frac{3C_4 M}{(n-k+\varepsilon^2)^{3/2}} L_\phi \varphi_{p-m,m}$$

and

$$|\mathbb{E}[(D^3 f_k(S_{k-p}) - D^3 f_k(S_{k-p-1})) * X_{k-m}] * \mathbb{E}[X_{k-j} \otimes X_k]|$$

$$\leq d^4 C \frac{3C_4 M^3}{(n-k+\varepsilon^2)^{3/2}} L_\phi \varphi_{p-m,0}.$$

The sum of these quantities over $(j, m, p)$ satisfying $1 \leq j \leq j+1 \leq m$ and $(r+2)m+1 \leq p \leq \ell$ is less than

$$\frac{\tilde{K}_4}{(n-k+\varepsilon^2)^{3/2}} L_\phi \left( \sum_{p=1}^{\lfloor \sqrt{n-k} \rfloor} p^2 \varphi_{p,0} + \sum_{p=1}^{\lfloor \sqrt{n-k} \rfloor} p \sum_{j=1}^{\lfloor p/(r+1) \rfloor} \varphi_{p,j} \right).$$

Thus we have

$$\sum_{j=1}^{\ell} \sum_{m=j+1}^{\ell} |\mathrm{Cov}(D^3 f_k(S_{k-(r+2)m-1})$$

$$(46) \qquad\qquad - D^3 f_k(S_{k-\ell-1}), X_{k-m} \otimes X_{k-j} \otimes X_k)| \mathbf{1}_{\{\ell \geq (r+2)m+1\}}$$

$$\leq 2 \frac{\tilde{K}_4 \alpha_{\sqrt{n-k}}}{(n-k+\varepsilon^2)^{3/2}} L_\phi$$

and

$$\sum_{j=1}^{\ell} \sum_{m=j+1}^{\ell} |\mathbb{E}[(D^3 f_k(S_{k-(r+2)m-1})$$

$$(47) \qquad\qquad - D^3 f_k(S_{k-\ell-1})) * X_{k-m}]$$

$$\qquad\qquad\qquad * \mathbb{E}[X_{k-j} \otimes X_k]| \mathbf{1}_{\{\ell \geq (r+2)m+1\}}$$

$$\leq 2 \frac{\tilde{K}_4 \alpha_{\sqrt{n-k}}}{(n-k+\varepsilon^2)^{3/2}} L_\phi.$$

(b) Let us now control the following quantity:

$$\sum_{j=1}^{\ell} \sum_{m=j+1}^{\ell} \mathbb{E}[D^3 f_k(S_{k-(r+2)m-1}) - D^3 f_k(S_{k-\ell-1})]$$

$$\qquad\qquad * \mathbb{E}[X_{k-m} \otimes X_{k-j} \otimes X_k] \mathbf{1}_{\{\ell \geq (r+2)m+1\}}.$$

If $\ell \geq (r+2)m + 1$ and $m \geq (r+2)j + 1$, we have

$$|\mathbb{E}[D^3 f_k(S_{k-(r+2)m-1})] * \mathbb{E}[(X_{k-m} - \mathbb{E}[X_{k-m}]) \otimes X_{k-j} \otimes X_k]|$$

$$\leq d^3 K_3 \left( \frac{A}{(n-k+\varepsilon^2)^{3/2}} + \frac{1}{n} \right) L_\phi 2CM \varphi_{m-j,j}.$$

If $\ell \geq (r+2)m + 1$ and $m \leq (r+2)j$, we have

$$|\mathbb{E}[D^3 f_k(S_{k-(r+2)m-1})] * \mathbb{E}[X_{k-m} \otimes X_{k-j} \otimes (X_k - \mathbb{E}[X_k])]|$$

$$\leq d^3 K_3 \left( \frac{A}{(n-k+\varepsilon^2)^{3/2}} + \frac{1}{n} \right) L_\phi C3M^2 \varphi_{j,0}.$$

Therefore, we have

$$\sum_{j=1}^{\ell} \sum_{m=j+1}^{\ell} |\mathbb{E}[D^3 f_k(S_{k-(r+2)m-1})] * \mathbb{E}[X_{k-m} \otimes X_{k-j} \otimes X_k] \mathbf{1}_{\{\ell \geq (r+2)m+1\}}|$$

$$\leq d^3 \tilde{K}_6 \left( \frac{A}{(n-k+\varepsilon^2)^{3/2}} + \frac{1}{n} \right)$$

(48)

$$\times L_\phi \left( \sum_{p=1}^{\lfloor \sqrt{n-k} \rfloor} \sum_{j=1}^{\lfloor p/(r+1) \rfloor} \varphi_{p,j} + \sum_{p=1}^{\lfloor \sqrt{n-k} \rfloor} p \varphi_{p,0} \right)$$

$$\leq d^3 \tilde{K}_6 \left( \frac{A \alpha_{\sqrt{n-k}}}{(n-k+\varepsilon^2)^{3/2}} + \frac{\beta_{\sqrt{n-k}}}{n} \right) L_\phi.$$

In the same way, we get

$$\sum_{j=1}^{\ell} \sum_{m=j+1}^{\ell} |\mathbb{E}[D^3 f_k(S_{k-\ell-1})] * \mathbb{E}[X_{k-m} \otimes X_{k-j} \otimes X_k] \mathbf{1}_{\{\ell \geq (r+2)m+1\}}|$$

(49)

$$\leq d^3 \tilde{K}_6 \left( \frac{A \alpha_{\sqrt{n-k}}}{(n-k+\varepsilon^2)^{3/2}} + \frac{\beta_{\sqrt{n-k}}}{n} \right) L_\phi.$$

This completes the proof of Lemma A.2.3. $\quad \square$

**A.3. End of the proof of Proposition A.1.** Let an integer $n \geq 9n_0$ and a real number $A \geq M$ be given. Let us suppose that property $(\mathcal{P}_n(A))$ is satisfied. Let a real number $\varepsilon \geq 1$ be given.

Let us recall that we have

$$\mathbb{E}[\phi(S_n + \varepsilon Y)] - \mathbb{E}[\phi(S_{n_0-1} + T_{n_0-1,n} + \varepsilon Y)] = \sum_{k=n_0}^{n} \Delta_k(f_{\phi,k,n,\varepsilon}),$$

with $T_{n_0-1,n} := \sum_{i=n_0}^{n} Y_i$. Then, according to Proposition A.1.3, we have

$$\sum_{k=n_0}^{n-\lfloor n/3 \rfloor -1} |\Delta_k(f_{\phi,k,n,\varepsilon})| \le K_1'' L_\phi,$$

where $K_1''$ only depends on $(d, C, M, r, (\varphi_{p,l})_{p,l})$.

On the other hand, according to Lemma A.2.2, we have

$$(50) \qquad \sum_{k=n-\lfloor n/3 \rfloor}^{n} |\Delta_{2,k}(f_{\phi,k,n,\varepsilon})| \le L_\phi K_2'' \left( 1 + A \sum_{l \ge 0} \frac{1}{(l+\varepsilon^2)^{3/2}} \right).$$

Moreover, according to Lemma A.2.3, for any integer $k = n - \lfloor \frac{n}{3} \rfloor, \ldots, n$, we have

$$|\Delta_{1,k}(f_{\phi,k,n,\varepsilon})| \le L_\phi \tilde{K} \bigg( \frac{A\alpha_{\sqrt{n-k}}}{(n-k+\varepsilon^2)^{3/2}} + \frac{\beta_{\sqrt{n-k}}}{n} + \frac{A\gamma_{\sqrt{n-k}}}{n-k+\varepsilon^2}$$
$$+ \frac{A\delta_{\sqrt{n-k}}}{\sqrt{n-k+\varepsilon^2}} + \varphi_{\lfloor \sqrt{n-k} \rfloor +1,0} \bigg),$$

with $\alpha_m = 1 + \sum_{p=1}^{\lfloor m \rfloor} p\zeta_p$, $\beta_m = 1 + \sum_{p=1}^{\lfloor m \rfloor} \zeta_p$, $\gamma_m = \sum_{p=\lceil m/(r+2)^2 \rceil}^{\lfloor m \rfloor} \zeta_p$ and $\delta_m = \sum_{p=\lceil m/(r+2) \rceil +1}^{+\infty} \frac{\zeta_p}{p}$ with $\zeta_p := p \max_{j=0,\ldots,\lfloor p/(r+1) \rfloor} \varphi_{p,j}$. We control independently each sum of these terms over $k \in \{n - \lfloor \frac{n}{3} \rfloor, \ldots, n\}$.

1. Control of the first term:

$$\sum_{k=n-\lfloor n/3 \rfloor}^{n} \frac{A\alpha_{\sqrt{n-k}}}{(n-k+\varepsilon^2)^{3/2}}$$

$$(51) \qquad = A \sum_{l=0}^{\lfloor n/3 \rfloor} \frac{1 + \sum_{p=1}^{\lfloor \sqrt{l} \rfloor} p\zeta_p}{(l+\varepsilon^2)^{3/2}}$$

$$\le A \left( \left( \sum_{l \ge 0} \frac{1}{(l+\varepsilon^2)^{3/2}} \right) + \sum_{p=1}^{\sqrt{n}} \left( \sum_{l \ge p^2} \frac{1}{(l+\varepsilon^2)^{3/2}} \right) p\zeta_p \right)$$

$$\le A \left( \left( \sum_{l \ge 0} \frac{1}{(l+\varepsilon^2)^{3/2}} \right) + \sum_{p \ge 1} \frac{2}{\sqrt{p^2+\varepsilon^2-1}} p\zeta_p \right).$$

2. Control of the second term:

$$(52) \qquad \sum_{k=n-\lfloor n/3 \rfloor}^{n} \frac{\beta_{\sqrt{n-k}}}{n} \le 1 + \sum_{p \ge 1} \zeta_p.$$

3. Control of the third term:

$$\sum_{k=n-\lfloor n/3 \rfloor}^{n} \frac{A\gamma_{\sqrt{n-k}}}{n-k+\varepsilon^2} \leq \sum_{l=0}^{\lfloor n/3 \rfloor} \frac{A}{l+\varepsilon^2} \sum_{p=\lceil \sqrt{l}/(r+2)^2 \rceil}^{\lfloor \sqrt{l} \rfloor} \zeta_p$$

$$\leq A \sum_{p=0}^{\lfloor \sqrt{n} \rfloor} \sum_{l=p^2}^{p^2(r+2)^4} \frac{1}{l+\varepsilon^2} \zeta_p.$$

Therefore, we have

$$(53) \qquad \sum_{k=n-\lfloor n/3 \rfloor}^{n} \frac{A\gamma_{\sqrt{n-k}}}{n-k+\varepsilon^2} \leq A \sum_{p\geq 0} \ln\left(\frac{p^2(r+2)^4+\varepsilon^2}{p^2-1+\varepsilon^2}\right) \zeta_p.$$

4. Control of the fourth term:

$$\sum_{k=n-\lfloor n/3 \rfloor}^{n} \frac{A\delta_{\sqrt{n-k}}}{\sqrt{n-k+\varepsilon^2}} \leq A \sum_{l=0}^{\lfloor n/3 \rfloor} \frac{1}{\sqrt{l+\varepsilon^2}} \sum_{p=\lfloor \sqrt{l}/(r+2) \rfloor + 1}^{+\infty} \frac{\zeta_p}{p}$$

$$\leq A \sum_{p\geq 1} \sum_{l=0}^{p^2(r+2)^2} \frac{1}{\sqrt{l+\varepsilon^2}} \frac{\zeta_p}{p}$$

$$\leq A \sum_{p\geq 1} 2\left(\sqrt{p^2(r+2)^2+\varepsilon^2} - \sqrt{\varepsilon^2-1}\right) \frac{\zeta_p}{p}.$$

Hence, we have

$$(54) \qquad \sum_{k=n-\lfloor n/3 \rfloor}^{n} \frac{A\delta_{\sqrt{n-k}}}{\sqrt{n-k+\varepsilon^2}} \leq A \sum_{p\geq 1} \frac{2(r+2)^2(1+p^2)}{\sqrt{p^2+\varepsilon^2}} \frac{\zeta_p}{p}.$$

5. Control of the fifth term:

$$\sum_{k=n-\lfloor n/3 \rfloor}^{n} \varphi_{\lfloor \sqrt{n-k} \rfloor + 1, 0} = \sum_{l=0}^{\lfloor n/3 \rfloor} \varphi_{\lfloor \sqrt{l} \rfloor + 1, 0}$$

$$\leq \sum_{p\geq 0} \#\{l : \lfloor \sqrt{l} \rfloor = p\} \varphi_{p+1, 0},$$

from which we get

$$(55) \qquad \sum_{k=n-\lfloor n/3 \rfloor}^{n} \varphi_{\lfloor \sqrt{n-k} \rfloor + 1, 0} \leq \sum_{p\geq 0} (2p+1) \varphi_{p+1, 0}.$$

**Acknowledgments.** I thank Jean-René Chazottes for having indicated the Kantorovich metric to me. I am grateful to Jean-Pierre Conze and to Yves Derriennic for our always interesting discussions. I also wish to thank the referees for their careful reading of this paper and for their pertinent remarks.

Département de Mathématiques
Faculté des Sciences et Techniques
Université de Bretagne Occidentale
29285 Brest Cedex
France
e-mail: francoise.pene@univ-brest.fr